\documentclass[12pt,a4paper]{article}

\usepackage{amsmath}    
\usepackage{amssymb}
\usepackage{graphicx}   
\usepackage{verbatim}   
\usepackage{color}      

\usepackage{hyperref}   

\usepackage{enumerate}
\usepackage{mathrsfs}
\usepackage{empheq}
\usepackage{bbold}

\usepackage[final]{showlabels}

\usepackage{amsthm}

\newtheorem{corollary}{Corollary}[section]
\newtheorem{theorem}{Theorem}[section]

\newcommand{\del}{\partial}
\renewcommand{\theta}{\vartheta}
\renewcommand{\phi}{\varphi}

\newcommand{\veccc}[3]{\left ( \begin{array}{c}#1\\#2\\#3\\ \end{array}\right )}

\newcommand{\dd}{\mathrm{d}}
\newcommand{\ee}{\mathrm{e}}

\renewcommand{\vec}{\mathbf}

\newcommand{\const}{\mathrm{const}}

\newcommand{\ii}{\mathbb{i}}

\renewcommand{\title}{On the Active Flux scheme for hyperbolic PDEs with source terms}

\newcommand{\authorFirst}{Wasilij Barsukow\footnote{Institute of Mathematics, Zurich University, 8057 Zurich, Switzerland}}
\newcommand{\authorTwo}{Jonas P.\ Berberich\footnote{\label{foot:wue}Wuerzburg University, Emil-Fischer-Strasse 40, 97074 Wuerzburg, Germany}}
\newcommand{\authorThree}{Christian Klingenberg\footnotemark[2]{}}

\begin{document}

\begin{center} \Large
\title

\vspace{1cm}

\date{}
\normalsize

\phantom{}

\authorFirst, \authorTwo, \authorThree
\end{center}

\begin{abstract}

The Active Flux scheme is a Finite Volume scheme with additional point values distributed along the cell boundary. It is third order accurate and does not require a Riemann solver: the continuous reconstruction serves as initial data for the evolution of the points values. The intercell flux is then obtained from the evolved values along the cell boundary by quadrature. This paper focuses on the conceptual extension of Active Flux to include source terms, and thus for simplicity assumes the homogeneous part of the equations to be linear. To a large part, the treatment of the source terms is independent of the choice of the homogeneous part of the system. Additionally, only systems are considered which admit characteristics (instead of characteristic cones). This is the case for scalar equations in any number of spatial dimensions and systems in one spatial dimension. Here, we succeed to extend the Active Flux method to include (possibly nonlinear) source terms while maintaining third order accuracy of the method. This requires a novel (approximate) operator for the evolution of point values and a modified update procedure of the cell average. For linear acoustics with gravity, it is shown how to achieve a well-balanced / stationarity preserving numerical method.



\end{abstract}

\section{Introduction}

Numerous phenomena of the physical world are modeled by hyperbolic balance laws (conservation laws augmented by source terms). This includes gas dynamics, the motion of water waves, plasma physics and even general relativity. Often physical modeling requires to include source terms, and conservation is modified due to creation or annihilation of some of the evolved quantities. Chemical reactions, for example, change the number density of a species and produce or absorb heat (i.e. internal energy). Gravity accelerates matter downwards and creates momentum. In the shallow water model describing the motion of a free water surface the bottom topography enters the equations through a source term. Rewriting the hydrodynamic equations in a different coordinate system (e.g. in polar coordinates) makes geometric source terms appear. 
All these applications require reliable numerical methods which are able to deal with source terms. 

Numerical methods for hyperbolic conservation laws with source terms first need to perform well in the homogeneous case. This means for example that they need to cope with discontinuities / weak solutions and with phenomena arising in multiple spatial dimensions, such as involutions and non-trivial stationary states. {This requirement has led \cite{eymann13,fan15} to suggest \emph{Active Flux}, an extension of the finite volume method. Additionally to the cell average, this scheme evolves point values located at the cell boundary. These are shared among neighbouring cells, which gives rise to a continuous reconstruction. The update of the point values is achieved by using an evolution operator that includes multi-dimensional information. The presence of the point values along the cell boundary then allows to compute the intercell flux via quadrature. Thus, Active Flux does not use Riemann solvers, while still evolving the cell average as one of the discrete degrees of freedom just as Finite Volume methods do. The additional (pointwise) degrees of freedom allow for the scheme to be of high order of accuracy on a compact stencil. It has been shown in \cite{barsukow18activeflux} that this scheme is stationarity preserving and vorticity preserving for linear acoustics without any fix. It is third order accurate. Extensions to nonlinear systems have been recently suggested e.g. in \cite{fan17,kerkmann18,barsukow19activeflux}. Active flux therefore seems to be promising for resolving many of the structure preservation problems that currently available methods are facing (an overview of existing methods for balance laws is given below).}

In view of the many applications that involve source terms, this paper therefore aims at deriving the necessary modifications for Active Flux to be applicable to balance laws while retaining its third order accuracy. Active flux for equations with a source term was considered in \cite{nishikawa16}, where for stationary problems the necessary quadratures could be chosen of lower order of accuracy (trapezoidal rule) than in the original Active Flux method from \cite{eymann13} (Simpson's rule) (see e.g. Eqn. (32) in \cite{nishikawa16}). For time-dependent problems, in \cite{nishikawa16} the reduced order of accuracy of these quadratures is remedied by using a high-order implicit time stepping method. The approach of the present work avoids sub-iterations and multi-step time integrators, and the high order in time is achieved through the choice of high order quadratures, that hardly entail any computational cost. Contrary to \cite{nishikawa16}, this paper presents a fully explicit method for hyperbolic problems with source terms that reverts to the original Active Flux scheme of \cite{eymann11} when the source term vanishes. As we aim at resolving the acoustic time scale, explicit time stepping is very efficient.

Including the source term requires a number of modifications. The homogeneous part of the equations therefore is for simplicity assumed to be a linear hyperbolic problem for which characteristics are available. This is the case for scalar equations in any number of spatial dimensions and for systems in one spatial dimension. For multi-dimensional systems, the concept of characteristics needs to be replaced by characteristics cones. In the homogeneous case, Active Flux has been used for this situation as well (\cite{eymann13,barsukow18activeflux}), but an extension to inhomogeneous systems in multi-d, and to nonlinear systems remains subject of future work. To a large part, the strategies presented in this paper will, however, remain valid when the homogeneous part of the equations is nonlinear as well, and even for nonlinear multi-dimensional systems.

As soon as a source term is added to a hyperbolic system, new stationary states arise which often are of particular interest. The stationarity is due to the flux divergence being equal to the source term. Many areas of application of balance laws involve studies of dynamics on top of such an equilibrium (e.g. astrophysics, meteorology, tsunami modeling, \ldots). This requires the numerical method to be very accurate on the stationary states in order to avoid spurious, artificial perturbations. Therefore the error of a numerical solution representing one of those stationary states should not increase with time, thus allowing the simulation to run for a long time {(see e.g. the review \cite{edelmann21})}. 

Numerical methods which achieve this are called \emph{well-balanced}, introduced in \cite{Greenberg1996}. They make sure that the discretization of the flux divergence and the discretization of the source term match, and that the numerical method keeps the desired stationary state exactly stationary for any resolution of the grid. The concept of well-balanced methods has been extensively used in the context of shallow water equations with non-flat bottom topography (e.g.\ \cite{Audusse2004,Bermudez1994,leveque98} and references therein). Here, the balance is the so-called lake-at-rest solution, which amounts to an algebraic condition and can thus be given explicitly. 

Another area in which well-balanced methods have high relevance is the simulation of hydrodynamic processes using compressible Euler equations with gravitational source term. The so-called hydrostatic state (stationary state with no velocity) is described by one PDE for two unknown functions. There are many hydrostatic states, depending on the additional thermodynamical relation that one chooses in order to close this PDE. The fact that the stationary state is itself given by a differential equation that cannot be immediately integrated makes well-balancing much more delicate in this context. There are two different ways which are currently used to construct well-balanced methods for the Euler equations with gravity. The first and more traditional way is to restrict the class of hydrostatic solutions which are balanced exactly or to choose a particular, but arbitrary hydrostatic state (e.g. \cite{Cargo1994,LeVeque2011,Desveaux2014,Chandra15,Berberich2016,Chertock2018a,Berberich2019,Berberich2019b}). This is advantageous in all those applications where the stationary state is known, and the evolution of perturbations around it shall be studied. If no information on the stationary state can be assumed, then the only way to proceed is to make sure that the stationary states of the numerical method are fulfilling some \emph{discretization} of the corresponding PDE (e.g.\ \cite{zenk14,Kaeppeli2016,Berberich2020}).

In this paper this latter approach is used. In the situation of the stationary states given by underdetermined PDEs, and not by algebraic equations, the relation between the discrete stationary states and the stationary states of the PDE has been studied in \cite{barsukow17a} for linear problems. It turns out that many standard numerical methods add diffusion even to those states that should remain stationary. The set of states that are actually kept stationary by such methods is very small (e.g. uniform constants). \emph{Stationarity preserving} methods, on the other hand, do not apply diffusion to discrete data which fulfill a discrete version of the PDE governing the stationary states. Stationarity preserving methods thus keep stationary a much larger set of initial data. Independently of how these discrete equations actually look like, it is their existence that makes a qualitative difference. In a non-stationarity-preserving method, initial data sampled from an analytic stationary state will decay due to the diffusion and become unrecognizable in the end. In a stationarity preserving method, these initial data will evolve towards one of the many discrete stationary states approximating the steady PDE, and will remain there forever (up to machine precision). The long-time numerical solution will then indeed approximate the analytic stationary state. For more details, see \cite{barsukow17a}. In this paper we understand the concept of well-balancing in this sense of stationarity preservation.

After extending the Active Flux scheme to include source terms, we construct a well-balanced Active Flux method for the equations of acoustics with gravity. The hydrostatic solutions of acoustics with gravity are comparable to those of the compressible Euler equations with gravity, since they are given via the same underdetermined differential equation. We show that the Active Flux scheme endowed with an exact evolution operator is intrinsically well-balanced. In practice, an approximate evolution operator needs to be used. Hence we introduce an approximate evolution operator which retains the well-balanced property. 

The paper is organized as follows: After the Active Flux scheme for homogeneous problems is introduced in section \ref{sec:af}, the modifications necessary for including source terms are discussed: Section \ref{sec:pointvalueswithsource} discusses the evolution operators necessary for the update of the point values and Section \ref{sec:updateaverage} is devoted to the modifications in the update of the average. Here, the focus lies on linear systems of equations with possibly nonlinear source terms in one spatial dimension and on linear advection in multiple spatial dimensions. Section \ref{sec:wellbal} discusses well-balancing of Active Flux for linear acoustics with gravity. Section \ref{sec:numeric} finally demonstrates numerically that the new method attains third order accuracy with linear and nonlinear source terms, can be used to compute Riemann problems, and displays well-balanced behavior for stationary states.

This work can be seen in the larger context of the quest for structure preserving numerical methods, of which well-balanced methods form an example. Extending these results to nonlinear hyperbolic equations with source terms and thus combining the structure preserving properties of Active Flux remains subject of future work. However, the procedures suggested in this paper are formulated with as little reference to the linearity of the equations as possible.

\section{The Active Flux scheme} \label{sec:af}

Consider the initial value problem for an $m \times m$ system of hyperbolic balance laws in $d$ spatial dimensions\footnote{In this paper, indices never denote derivatives. Boldface symbols denote vectors that have the same dimension as the space.}
\begin{align}
 \del_t q + \nabla \cdot \vec f(q) =&\; s(q) & q &:\mathbb R^+_0 \times \mathbb R^d \to \mathbb R^m \label{eq:conslaw},\;
 f,s: \mathbb R^m \to \mathbb R^m \\ 
 q(0, \vec x) :=&\; q_0(\vec x)
\end{align}
This section reviews the general idea of the Active Flux scheme. {Instead of introducing jumps at every cell interface, as is customary for finite volume schemes, Active Flux employs a continuous reconstruction and evolves point values at the cell interfaces independently. These point values are shared by the adjacent cells. Thus, despite evolving a cell average Active Flux does not require a numerical flux function, as there is no Riemann Problem to solve. Given the point values, the update of the cell average is immediately possible by performing flux quadrature in time and along the cell interface. The distribution of degrees of freedom is discussed in section \ref{ssec:dofs}, and the update of the average in section \ref{ssec:reviewfinvol}. What remains, is the update of the point values. To this end, an IVP is solved (approximately) with the initial data given by the globally continuous reconstruction. This is very different from the usual approach of finite volume schemes and is described in sections \ref{ssec:reviewevo}--\ref{ssec:reviewrecon}.} Some of the details of the (approximate) evolution operator then depend on the particular equation that is to be solved. After the general concept is outlined, the details that make it applicable to hyperbolic balance laws are discussed in sections \ref{sec:pointvalueswithsource} and \ref{sec:updateaverage}.

\subsection{Degrees of freedom in the Active Flux scheme} \label{ssec:dofs}

\newcommand{\node}{}
\newcommand{\eh}{}
\newcommand{\ev}{}

\begin{figure}
 \centering
 \includegraphics[width=0.6\textwidth]{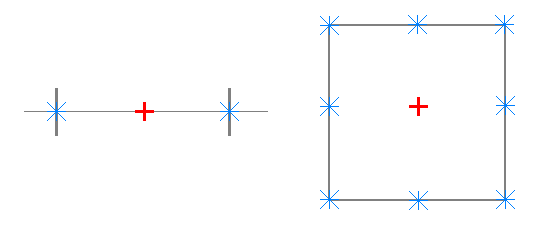}
  \caption{The degrees of freedom used for Active Flux. Stars indicate the location of point values, and the cross (placed in the center symbolically) refers to the cell average. \emph{Left}: One spatial dimension. \emph{Right}: Two spatial dimensions.}
 \label{fig:dofactivelux}
\end{figure}

The Active Flux scheme (\cite{eymann13,barsukow18activeflux}, first introduced in \cite{vanleer77}) is an extension of the finite volume scheme. The Active Flux scheme evolves both the cell average and point values which are distributed along the cell boundary. In particular, here the following two choices are considered (see Figure \ref{fig:dofactivelux}):
\begin{itemize}
 \item In one spatial dimension, there is a point value $q_{i+\frac12}$ located at each cell interface $x_{i+\frac12}$. Thus every cell has access to one cell average $\bar q_i$ and two point values at its interfaces. 
 \item On Cartesian grids in two spatial dimensions, there is a point value $q_{i+\frac12,j}$, $q_{i,j+\frac12}$ at each edge midpoint and one at each node $q_{i+\frac12,j+\frac12}$. Every cell has access to one cell average $\bar q_{ij}$ and 8 point values distributed along the cell interface.
\end{itemize}

{Note that the point values at cell interfaces are shared by the adjacent cells. As will be seen in the following, the reconstruction is globally continuous and no Riemann Problems arise.} In one spatial dimension, on average there are 2 degrees of freedom per cell: 1 cell average and 2 interface values shared each by 2 cells. In two spatial dimensions, in the setup described above, there are 4 degrees of freedom per cell: 1 cell average, 4 edge values, each shared by two cells and 4 node values each shared by 4 cells.

Note also that Active Flux does not use a staggered grid. The degrees of freedom at the cell boundaries are not averages over staggered volumes, but point values. This also explains why there is no notion of a conservative update for these, because this concept only applies to averages. The update of the cell average in the Active Flux method is, of course, conservative (see below).

\subsection{Update of the cell average} \label{ssec:reviewfinvol}

As the Active Flux scheme is an extension of the finite volume scheme, given a {time-step-average of the flux through the cell interface}, the update of the average happens in the same way as for finite volume schemes. {As there is a point value located at the cell interface, a Riemann Solver is not required to obtain the flux.} In this section, this finite volume aspect of Active Flux is described in an arbitrary number of spatial dimensions. 

Consider the computational domain to be subdivided into polygonal computational cells. Upon integration of \eqref{eq:conslaw} over one time step $[t^n, t^n + \Delta t]$ and over one computational cell $\mathcal C$ one obtains an evolution equation for the cell average $\bar q_\mathcal C := \frac{1}{|\mathcal C|} \int_{\mathcal C} \dd \vec x \, q(t, \vec x)$:
\begin{align*}
  \frac{\bar q^{n+1}_{\mathcal C} - \bar q^n_{\mathcal C}}{\Delta t}  + \frac{1}{|\mathcal C|} \frac{1}{\Delta t} \int\limits_{t^n}^{t^n + \Delta t} \!\!\!\! \dd t \, \int_{\del \mathcal C} \dd \sigma \, \vec n \cdot \vec f(q(t, \vec x)) &= \\
  \frac{1}{\Delta t} \int\limits_{t^n}^{t^n + \Delta t} \!\!\!\! \dd t \frac{1}{|\mathcal C|} \,& \int_{\mathcal C} \dd \vec x \, s(q(t , \vec x)) 
\end{align*}
Here, as usual, the index of the time step is placed as a superscript and $q^n_{\mathcal C}$ denotes the average in cell $\mathcal C$ at time $t^n$. The boundary $\del \mathcal C$ consists of edges $e$, such that one can rewrite
\begin{align*}
  \frac{\bar q^{n+1}_{\mathcal C} - \bar q^n_{\mathcal C}}{\Delta t}  + \frac{1}{|\mathcal C|} \frac{1}{\Delta t} \int\limits_{t^n}^{t^n + \Delta t} \!\!\!\! \dd t \, \sum_{e \subset \del \mathcal C} \int_{e} \dd \sigma \, \vec n_e \cdot \vec f(q(t, \vec x)) &=  \\ \frac{1}{\Delta t} \int\limits_{t^n}^{t^n + \Delta t} \!\!\!\! \dd t \frac{1}{|\mathcal C|} \, &\int_{\mathcal C} \dd \vec x \, s(q(t , \vec x)) 
\end{align*}
The vector $\vec n_e$ is the outward unit normal of edge $e$. This expression, so far exact, becomes a finite volume scheme upon replacing the exact normal flux and source averages by suitable approximations $\hat f_e$ and $\hat s_\mathcal C$:
\begin{align}
  \frac{\bar q^{n+1}_{\mathcal C} - \bar q^n_{\mathcal C}}{\Delta t} & + \frac{1}{|\mathcal C|}  \sum_{e \subset \del \mathcal C} |e| \hat f_e = \hat s_\mathcal C  \label{eq:scheme}\\
  \intertext{with}
  \hat f_e &\simeq \frac{1}{\Delta t} \int\limits_{t^n}^{t^n + \Delta t} \!\!\!\! \dd t \, \frac{1}{|e|} \int_{e} \dd \sigma \, \vec n_e \cdot \vec f(q(t, \vec x)) \label{eq:numflux}\\
  \hat s_\mathcal C &\simeq \frac{1}{\Delta t} \int\limits_{t^n}^{t^n + \Delta t} \!\!\!\! \dd t \, \frac{1}{|\mathcal C|} \, \int_{\mathcal C} \dd \vec x \, s(q(t , \vec x)) \label{eq:numsource}
\end{align}

Usual finite volume schemes introduce a (piecewise continuous) reconstruction of the averages, and obtain the numerical flux by an exact or approximate short-time evolution of this reconstruction. For example, introducing a piecewise constant function whose averages match the given cell averages, and solving the Riemann problems at the cell interfaces allows to compute a numerical flux.

The Active Flux scheme does not need this. Indeed, the point values along the boundary can be used to immediately approximate \eqref{eq:numflux}--\eqref{eq:numsource} by quadrature. The desired properties (most importantly the desired order of accuracy) of the resulting scheme dictate the number of point values along each edge and also the points in time at which these point values need to be available. 

The source term also contributes to the update of the cell average. The quadrature necessary to approximate the source term average \eqref{eq:numsource} to sufficient order in space and time is suggested in this paper for the first time and discussed in section \ref{sec:updateaverage}.

\subsection{Update of the point values} \label{ssec:reviewevo}

The cell average update, and in particular the computation of the intercell fluxes, requires accurate point values at the cell boundary to be available.

First consider the case where the source term vanishes: $s = 0$. For third order of accuracy, the integrals in \eqref{eq:numflux} need to be approximated by Simpson's rule. For the integration in space this can easily be achieved using the available point values at each cell interface as described in section \ref{ssec:dofs}. For the integration in time all point values need to be available at $t^n, t^n + \frac{\Delta t}{2}$ and $t^n + \Delta t$. Altogether this yields a space-time Simpson rule. 

In order to obtain sufficiently accurate time evolved point values, in \cite{vanleer77} it has been suggested to reconstruct the data and to use an exact evolution operator. An exact evolution operator generally is unavailable for nonlinear problems, and therefore in \cite{fan17,kerkmann18,barsukow19activeflux} approximate evolution operators have been proposed. Even for linear systems of hyperbolic balance laws it is generally very difficult to obtain closed-form exact evolution operators, as is shown in section \ref{ssec:evoacgrav}. Therefore the point values in the Active Flux scheme shall be evolved using a sufficiently high order \emph{approximate} evolution operator applied to a reconstruction of the discrete data. An exact evolution operator provides the necessary upwinding in order to guarantee stability, and an approximate evolution operator needs to do the same. The approximate evolution operator is introduced in section \ref{ssec:evolinhyp}.

\subsection{Reconstruction} \label{ssec:reviewrecon}

The reconstruction shall interpolate the point values and its average over the computational cell shall match the given cell average. In the following, to simplify notation, in one spatial dimension a uniform grid is assumed, although the reconstruction can immediately be generalized to nonuniform grids. In two spatial dimensions, a Cartesian grid is used. See \cite{eymann13} for a reconstruction on triangular grids. As mentioned in section \ref{ssec:dofs}, in one spatial dimension every cell has access to 3 degrees of freedom which makes a parabolic reconstruction natural. With the above-mentioned setup it is unique and reads (\cite{vanleer77,fan15})
\begin{align}
 q_{\text{recon},i}(x) &= -3 (2 \bar q_i - q_{i-\frac12} - q_{i+\frac12}) \frac{(x-x_i)^2}{\Delta x^2} \label{eq:parabolicrecon} \\\nonumber &+ (q_{i+\frac12} - q_{i-\frac12}) \frac{x-x_i}{\Delta x} + \frac{6 \bar q_i - q_{i-\frac12} - q_{i+\frac12}}{4}  \qquad x \in [x_{i-\frac12},x_{i+\frac12}]
\end{align}
In two spatial dimensions, in the setup as described above, every cell has access to 9 degrees of freedom, and there is a unique biparabolic reconstruction, which reads
\begin{align}
\begin{split}
 q_{\text{recon},ij}&(\xi \Delta x, \eta \Delta y) := \frac{9}{4} \bar q_{ij} \left(-1+4 \xi^2\right) \left(-1+4 \eta^2\right)\\
 &-\frac{1}{4} {q_\text W} \left(-1-4 \xi+12 \xi^2\right) \left(-1+4 \eta^2\right)\\
 &-\frac{1}{4} {q_\text E} \left(-1+4 \xi+12 \xi^2\right) \left(-1+4 \eta^2\right)\\
 &-\frac{1}{4} {q_\text S} \left(-1+4 \xi^2\right) \left(-1-4 \eta+12 \eta^2\right)\\
 &-\frac{1}{4} {q_\text N} \left(-1+4 \xi^2\right) \left(-1+4 \eta+12 \eta^2\right)\\
 &+\frac{1}{16} {q_\text{SW}} (-1+2 \xi) (-1+2 \eta) (-1-2 \eta +2 \xi (-1+6 \eta))\\
 &+\frac{1}{16} {q_\text{SE}} (1+2 \xi) (-1+2 \eta) (1+2 \eta+2 \xi (-1+6 \eta))\\
 &+\frac{1}{16} {q_\text{NW}} (-1+2 \xi) (1+2 \eta) (1-2 \eta+2 \xi (1+6 \eta))\\
 &+\frac{1}{16} {q_\text{NE}} (1+2 \xi) (1+2 \eta) (-1+2 \eta+2 \xi (1+6 \eta))
 \end{split} \label{eq:recon}
\end{align}
with $\xi := \frac{x-x_i}{\Delta x}$, $\eta := \frac{y-y_j}{\Delta y}$, $\xi \in \left[ -\frac12, \frac12 \right ]$, $\eta \in \left[ -\frac12, \frac12 \right ]$ and
\begin{align}
 q_\text{NE} &= q\node_{i+\frac12,j+\frac12} & q_\text{NW} &= q\node_{i-\frac12,j+\frac12} & q_\text{SW} &= q\node_{i-\frac12,j-\frac12} & q_\text{SE} &= q\node_{i+\frac12,j-\frac12} \\ 
 q_\text{N} &= q\eh_{i,j+\frac12} & q_\text{S} &= q\eh_{i,j-\frac12} & q_\text{E} &= q\ev_{i+\frac12,j} & q_\text{W} &= q\eh_{i-\frac12,j} 
\end{align}

Note that both reconstructions are globally continuous {and no Riemann Problems are introduced}. The reconstruction, however, is generally not continuously differentiable at the cell interfaces.

\subsection{Overview of the algorithm}

The overall algorithm of Active Flux is as follows:

\begin{enumerate}
 \item Given cell averages and point values, compute a reconstruction according to section \ref{ssec:reviewrecon}.
 \item Use the reconstruction as initial data in the update of the point values. The choices of evolution operators considered so far are discussed in section \ref{ssec:reviewevo} and evolution operators in presence of source terms are suggested in section \ref{ssec:evolinhyp} below.
 \item Given the updated point values along the cell interfaces, compute the intercell fluxes via quadrature (sections \ref{ssec:reviewfinvol} and \ref{sec:updateaverage} for the homogeneous and the inhomogeneous cases, respectively).
 \item Update the cell averages via \eqref{eq:scheme}.
\end{enumerate}

{The computations performed in the Active Flux algorithm are similar in structure and amount to high order Finite Volume methods, leading to similar time consumption in practice. The latter require a repeated evaluation of the reconstruction and of the numerical flux function for the individual steps of a time integrator (e.g. a Runge-Kutta method), while Active Flux performs several evaluations of the evolution operator to compute values for the flux quadrature in time (without recomputing the reconstruction). The shared degrees of freedom lead to lower memory usage in comparison to e.g. Discontinuous Galerkin (DG) methods.}

A CFL-type condition arises in the update of the point values: the domain of dependence of the evolution operator needs to be contained in the neighboring cells. Denoting by $\lambda_\text{max}$ the maximum speed of propagation, the time step needs to be chosen as
\begin{align}
 \Delta t \leq \frac{L_\text{min}}{\lambda_\text{max}}
\end{align}
where $L_\text{min} = \Delta x$ in one spatial dimension, and $L_\text{min} =\frac12 \min(\Delta x, \Delta y)$ in two spatial dimensions, if the point values are distributed as described in section \ref{ssec:dofs}. We introduce the CFL number as $\Delta t \lambda_\text{max} / L_\text{min}$.

\section{Evolution of the point values in presence of a source term} \label{sec:pointvalueswithsource}

The evolution of the point values needs to account for the source term. Additionally, in this paper a special focus shall lie on structure preservation properties of the resulting scheme. In the homogeneous case such properties have been observed upon usage of an exact evolution operator (\cite{barsukow18activeflux}). In presence of a source term, one needs to use an approximate evolution operator (section \ref{ssec:evolinhyp}), but should nevertheless aim at making it such that it does not spoil structure preservation (see section \ref{sec:wellbal}).

For certain equations, the inhomogeneous problem admits an exact solution (sections \ref{ssec:evolinadv}--\ref{ssec:evoacgrav}). This is valuable in order to assess specific properties of the numerical method later.

\subsection{Linear advection with a source term in multiple spatial dimensions} \label{ssec:evolinadv}

Consider a scalar equation ($m=1$) and $\vec f(q) = \vec U q$ with $\vec U \in \mathbb R^d$. Then
\begin{align}
 \del_t q + \vec U \cdot \nabla  q &= s(q)
\end{align}
amounts to the ODE 
\begin{align}
 \frac{\dd}{\dd t} q &= s(q)
\end{align}
along the straight characteristic of velocity $\vec U$. This ODE can be easily solved analytically:
\begin{align}
 \int_{q_0(\vec x - \vec U t)}^{q(t,\vec x)} \frac{\dd p}{s(p)} &= t
\end{align}
E.g. for $s(q) = \kappa q$ this yields $\ln \frac{q(t,\vec x)}{q_0(\vec x - \vec U t)} = \kappa t$, or
\begin{align}
 q(t,\vec x) &= q_0(\vec x - \vec U t) \exp(\kappa t) \label{eq:linadvlinsourceexact}
\end{align}
and for $s(q) = \kappa q^B$, $B \neq 1$
\begin{align}
 q(t,\vec x)  &= \Big ( (q_0(\vec x - \vec U t))^{1-B} + (1 - B) \kappa t \Big )^{\frac{1}{1-B}} \label{eq:linadvnonlinsourceexact}
\end{align}

\subsection{Linear acoustics with gravity in one spatial dimension}
\label{ssec:evoacgrav}

This section has threefold purpose. First, it introduces the acoustic equations with a gravity source term, which form a very useful system for the study of structure preservation of numerical methods. This is the set of equations for which a well-balanced method is derived in \ref{sec:wellbal}. This section also demonstrates the difficulties of finding an exact solution to an inhomogeneous system even if it is linear. Finally, the exact solution derived here is used later in order to assess the accuracy of the numerical method.

The equations of linear acoustics in one spatial dimension endowed with a gravity source term read:
\begin{align}
 \del_t \rho + \del_x v &= 0 \label{eq:acgrav1}\\
 \del_t v + \del_x p &= \rho g \qquad g \in \mathbb R\\
 \del_t p + c^2 \del_x v &= 0 \label{eq:acgrav3}
\end{align}

The corresponding homogeneous problem (linear acoustics) is the linearization of the Euler equations around the background state of constant density $\rho_\text{bg}= 1$, constant pressure $p_\text{bg}$ and vanishing velocity. Then the speed of sound $c = \sqrt{\frac{\gamma p_\text{bg}}{\rho_\text{bg}}}$ is a constant ($\mathbb R \ni \gamma > 1$). The full system \eqref{eq:acgrav1}--\eqref{eq:acgrav3} can be understood as a particular kind of a linearization of the Euler equations with gravity\footnote{Note that often the energy equation is written with a source term $\rho g v$. This source term is unnecessary, as it can be removed by redefining the notion of total energy. When the total energy includes the potential energy $-\rho g x$ due to gravity, the conservation form of the energy equation is restored. The source term in the momentum equation remains.}
\begin{align}
 \del_t \rho + \del_x (\rho v) &= 0 \label{eq:eulergrav1}\\
 \del_t (\rho v) + \del_x (\rho v^2 + p) &= \rho g \label{eq:eulergrav2}\\
 \del_t e + \del_x (v (e  + p  ) ) &= 0 \label{eq:eulergrav3}\\
 e = \frac{p}{\gamma-1} + \frac12 \rho v^2 &-  \rho g x
\end{align}

The static (stationary and $v=0$) states of \eqref{eq:eulergrav1}--\eqref{eq:eulergrav3} are governed by $\del_x p = \rho g$. This equation can only be solved if e.g. $\rho$ is given as a function of $x$, or if another relation is provided between any two of the variables $p, \rho, e$. This multitude of possible stationary states is reflected in the linearization \eqref{eq:acgrav1}--\eqref{eq:acgrav3}. (This is the reason for this particular choice of a linearization.) Observe that stationary states of \eqref{eq:acgrav1}--\eqref{eq:acgrav3} also are governed by $\del_x p = \rho g$ and that $p$ can only be computed if $\rho$ is given as a function of $x$, or if an additional relation is provided that links $\rho$ and $p$. This is an example of a so-called non-trivial stationary state as introduced in \cite{barsukow17a}. Examples of stationarity preserving schemes for \eqref{eq:acgrav1}--\eqref{eq:acgrav3} have been discussed in \cite{barsukow18thesis}.

The exact solution of \eqref{eq:acgrav1}--\eqref{eq:acgrav3} is studied in the Appendix \ref{sec:excatsolutionac}. This solution is not part of the suggested method but only serves auxiliary purposes, such as accuracy checks. However it illustrates the difficulties encountered when solving linear systems with sources. To the authors' knowledge the exact solution to \eqref{eq:acgrav1}--\eqref{eq:acgrav3} is not available in the literature so far.

\subsection{Runge-Kutta method for linear systems with a source} \label{ssec:evolinhyp}

Consider an $m \times m$ linear system in characteristic variables:
\begin{align}
 (\del_t + \lambda_\ell \del_x) Q_\ell &= S_\ell(Q_1, \ldots, Q_m) \qquad \ell = 1, \ldots, m\\
 \nonumber Q_\ell \colon \mathbb R^+_0 \times \mathbb R &\to \mathbb R \qquad \lambda_\ell \in \mathbb R \qquad S_\ell \colon \mathbb R^m \to \mathbb R
\end{align}
From now on, the capital letter $Q$ denotes the characteristic variables of this particular system, whereas $q$ continues to denote a generic variable.

Recall the following theorem from \cite{barsukow19activeflux}:
\begin{theorem}\label{thm:order}
 Assume a hyperbolic CFL condition $\Delta t / \Delta x \to \const$ as $\Delta t \to 0$. If the approximate evolution $Q^{\rm approx}(t, x)$ approximates the exact solution $Q(t, x)$ for fixed $x$ at least as 
\begin{align}
 Q^{\rm approx}(t, x) = Q(t, x) + \mathcal O(t^3)
\end{align}
and the quadrature rules used to approximate \eqref{eq:numflux}--\eqref{eq:numsource} yield the exact value up to an error of $\mathcal O(\Delta t^\alpha \Delta x^\beta)$, $\alpha + \beta \geq 3$ then Active Flux formally achieves third order accuracy.
\end{theorem}

Note that the simple approach of evolving each component of the source term along its associated characteristic
\begin{align}
Q_\ell(t, x) &\simeq Q_{\ell,0}(x - \lambda_\ell t) + t S_\ell(Q_{1,0}(x - \lambda_\ell t), \ldots, Q_{m,0}(x - \lambda_\ell t)) \qquad \ell = 1, \ldots, m
\end{align}
fails to be accurate enough (the error is $\mathcal O(t^2)$ instead of $\mathcal O(t^3)$).

Recall the second order Runge-Kutta method for the ordinary differential equation
\begin{align}
 \frac{\dd}{\dd t} q(t) &= s(t, q(t)) & q &: \mathbb R^+_0 \to \mathbb R  
\end{align}
\begin{align}
 q^{(1)}(\alpha t) &= q(0) + \alpha t s(0, q(0))\\
 q(t) &= q(0) + t \left( 1 - \frac{1}{2\alpha} \right ) s(0, q(0)) + t \frac{1}{2\alpha} s(\alpha t, q^{(1)}(\alpha t)) + \mathcal O(t^3)
\end{align}
for any $\alpha \in (0,1)$. In particular choosing $\alpha = \frac12$ (midpoint method) involves a predictor value at half time step. This can be taken as inspiration for constructing a sufficiently accurate approximate evolution operator:

\begin{figure}
 \centering
 \includegraphics[width=0.8\textwidth]{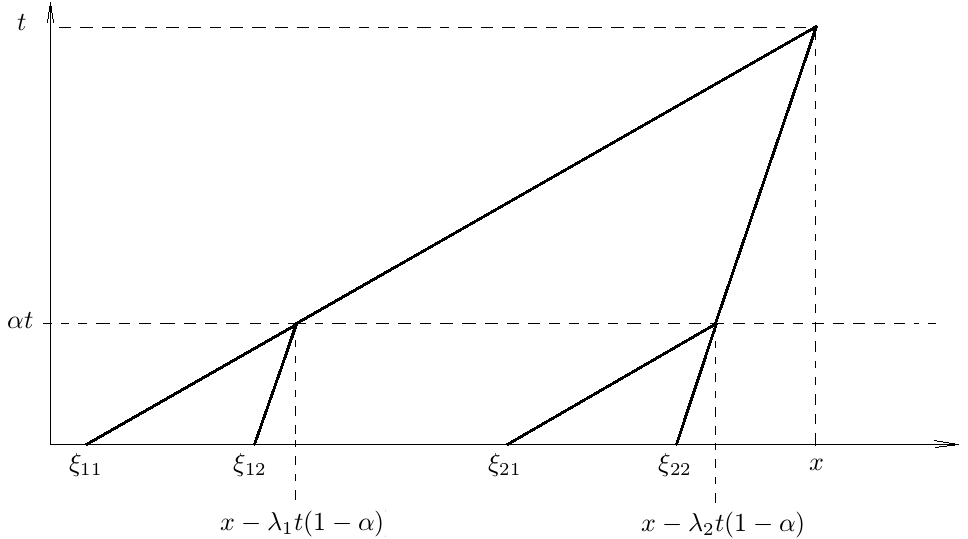}
 \caption{Illustration of the intermediate solutions and the involved characteristics for the first step in the Runge-Kutta scheme.}
 \label{fig:rk2sketch}
\end{figure}

\begin{theorem}[RK2 evolution operator] \label{thm:rk2}
Choose (see Figure \ref{fig:rk2sketch})
\begin{align}
 \xi_{\ell k} &:= x - \lambda_\ell  t(1-\alpha) - \lambda_k \alpha t\\
 Q_{k\ell }^* &:= Q_{k,0}(\xi_{\ell k}) + \alpha t S_k(Q_{1,0}(\xi_{\ell k}), \ldots, Q_{m,0}(\xi_{\ell k})) \qquad k,\ell  = 1, \ldots, m
\end{align}
and
\begin{align}
 Q_\ell^{(1)}(t, x) := Q_{\ell ,0}(x-\lambda_\ell  t) &+ \left( 1 - \frac{1}{2\alpha} \right ) S_\ell (Q_{1,0}(x-\lambda_\ell  t), \ldots, Q_{m,0}(x-\lambda_\ell  t)) t \\&+ \frac{t}{2\alpha} S_\ell \Big(Q_{1\ell }^*, \ldots, Q_{m\ell }^*\Big) \qquad \ell  = 1, \ldots, m \label{eq:rk2solution}
\end{align}

Then, for all $\alpha \in (0,1)$
\begin{align}
Q_\ell^{(1)}(t, x) &= Q_\ell(t, x) + \mathcal O(t^3) \qquad \ell = 1, \ldots ,m
\end{align}
\end{theorem}
Note that $Q_{\ell j}^*$ approximates $Q_\ell(\alpha t, x - \lambda_j t (1-\alpha))$.
\begin{proof}
By explicitly computing the first three terms of the Taylor series in $t$ one confirms the statement. The exact solution is
\begin{align}
 Q_\ell(t, x) &= Q_{\ell,0}(x) + t \del_t Q_\ell\Big |_{t=0}+ \frac{t^2}{2} \del_t^2 Q_\ell\Big |_{t=0} + \mathcal O(t^3)\\
 &=Q_{\ell,0}(x) + t (S_{\ell,0} - \lambda_\ell \del_x Q_{\ell,0}) \\\nonumber &+ \frac{t^2}{2} \left( \sum_k \frac{\del S_\ell}{\del Q_k} \Big (S_{k,0} - (\lambda_k + \lambda_\ell) \del_x Q_{k,0}\Big ) + \lambda_\ell^2 \del_x^2 Q_{\ell,0}) \right)  + \mathcal O(t^3)
\end{align}
where $S_{\ell,0}$ denotes
\begin{align}
 S_{\ell,0} := S_\ell(Q_{1,0}(x), \ldots,Q_{m,0}(x))
\end{align}
and $\frac{\del S_\ell}{\del Q_k}$ also is evaluated at $x$. Note that it has been used that $\del_x \lambda_\ell = 0$ (i.e. that the homogeneous system is linear), but the source $S$ can be any differentiable function of $Q$. 

Expand now \eqref{eq:rk2solution} ($\ell = 1, \ldots, m$):
\begin{align}
 \del_t Q_{k\ell}^*\Big |_{t=0} &= -(\lambda_\ell (1-\alpha) + \lambda_k \alpha) \del_x Q_{k,0} + \alpha  S_{k,0} \\
 \del_t Q_\ell^{(1)}(t, x) &= -\lambda_\ell \del_x Q_{\ell,0}(x-\lambda_\ell t) \\
 &+ \left( 1 - \frac{1}{2\alpha} \right ) \left( t\sum_k \frac{\del S_\ell}{\del Q_k} \del_x Q_{k,0}(x-\lambda_\ell t) (-\lambda_\ell) \right . \\& \left . \phantom{mmm\frac{\del S_\ell}{\del Q_k}}  + S_\ell(Q_{1,0}(x-\lambda_\ell t), \ldots, Q_{m,0}(x-\lambda_\ell t)) \right ) \\
 &+ \frac{1}{2\alpha} \left( t \sum_k \frac{\del S_\ell}{\del Q_k} \del_t Q_{k\ell}^*  + S_\ell\Big(Q_{1\ell}^*, \ldots, Q_{m\ell}^*\Big) \right )\\
 &\overset{t=0}{=} -\lambda_\ell \del_x Q_{\ell,0} +  S_{\ell,0}\\
 \del_t^2 Q_\ell^{(1)}(t, x) \Big |_{t=0} &= \lambda_\ell^2 \del_x^2 Q_{\ell,0} + \left( 1 - \frac{1}{2\alpha} \right ) \left(2 \sum_k \frac{\del S_\ell}{\del Q_k} \del_x Q_{k,0} (-\lambda_\ell) \right ) \\
 &+ \frac{1}{2\alpha} \left( 2 \sum_k \frac{\del S_\ell}{\del Q_k} \del_t Q_{k\ell}^* \Big |_{t=0}  \right )\\
 &= \lambda_\ell^2 \del_x^2 Q_{\ell,0} 
  - \sum_k \frac{\del S_\ell}{\del Q_k} \Big (  \del_x Q_{k,0} \left( \lambda_\ell  + \lambda_k \right ) -  S_{k,0} \Big ) 
\end{align}
Obviously the two Taylor series agree up to terms $\mathcal O(t^3)$, which proves the statement.

\end{proof}

\begin{corollary}[Midpoint method]
If $\alpha = \frac12$, then for $\ell,k = 1, \ldots m$

\begin{align}
 \xi_{\ell j} &:= x - (\lambda_\ell + \lambda_j ) \frac{t}{2}\\
 Q_{k\ell}^* &:= Q_{k,0}(\xi_{\ell k}) + \frac{t}{2} S_k(Q_{1,0}(\xi_{k\ell}), \ldots, Q_{m,0}(\xi_{k\ell}))\\
 Q_\ell^{(1)}(t, x) &:= Q_{\ell,0}(x-\lambda_\ell t) + t S_\ell\Big(Q_{1 \ell}^*, \ldots, Q_{m \ell}^*\Big) 
\end{align}

\end{corollary}

\begin{corollary}[RK2 evolution operator for a scalar equation] \label{cor:scalarrk2}
For a scalar equation
\begin{align}
 (\del_t + \lambda \del_x) Q &= S(Q)
\end{align}
the algorithm reads
\begin{align}
 \xi &:= x - \lambda t 
\end{align}
and
\begin{align}
 Q^{(1)}(t, x) &:= Q_{0}(x-\lambda t) 
 + \left( 1 - \frac{1}{2\alpha} \right ) S( Q_{0}(x-\lambda t)) t \\
 &+ \frac{t}{2\alpha} S\Big(  Q_{0}(\xi) + \alpha t S(Q_{0}(\xi)) \Big)
\end{align}
\end{corollary}

For the equations \eqref{eq:acgrav1}--\eqref{eq:acgrav3} of linear acoustics with gravity, $\lambda_1 = c = -\lambda_2, \lambda_3 = 0$. The characteristic variables are
\begin{align}
 Q_1 &= \frac{p + c v}{2} & Q_2 &= \frac{p-cv}{2} & Q_3 = -\frac{p}{c^2} + \rho  \label{eq:acgravcharvar}
\end{align}
and the gravity source term then is 
\begin{align}
 S_1 = -S_2 &= \frac{g}{2c} (Q_1 + Q_2) +  \frac{cg}{2} Q_3  & S_3 &= 0
\end{align}

\section{Update of the cell average in presence of a source term} \label{sec:updateaverage}

The update of the cell average needs to include the space-time average of the source term according to \eqref{eq:scheme} {of section \ref{ssec:reviewfinvol}}. This space-time average needs to be approximated by a suitable quadrature / approximation with sufficient order of accuracy. Active flux has a strong focus on providing discrete degrees of freedom along the boundary which allow to perform a quadrature along the boundary. However, the evaluation of the source term for the update of the cell average involves an averaging over the cell volume. It is more difficult to achieve the desired order of accuracy here, as the setup lacks the quadrature points that would have been natural for this task. A quadrature formula adapted to the geometry of the Active Flux method is derived here.

\subsection{One spatial dimension}

The approximation \eqref{eq:numsource}
\begin{align}
 \hat s_\mathcal C &\simeq \frac{1}{\Delta t} \int\limits_{t^n}^{t^n + \Delta t} \!\!\!\! \dd t \, \frac{1}{|\mathcal C|} \, \int_{\mathcal C} \dd \vec x \, s(q(t , \vec x))
\end{align}
of the source term in \eqref{eq:scheme} requires a space-time quadrature that is exact for parabolic functions. The natural candidate would be Simpson's rule in both space and time (as used for the numerical flux), but there are not enough quadrature points for it. For example in one spatial dimension, the available information is

\setlength{\tabcolsep}{8.8pt}
\renewcommand{\arraystretch}{2.0}
\begin{center}
\begin{tabular}{c|ccc}
 $t^{n+1}$ &$q_{i-\frac12}^{n+1}$ & &$q_{i+\frac12}^{n+1}$\\
 $t^{n+\frac12}$ &$q_{i-\frac12}^{n+\frac12} $   && $q_{i+\frac12}^{n+\frac12}$\\
$t^{n}$ & $q_{i-\frac12}^{n+1}$ & $\boxed{\bar q_i^n}$   & $q_{i+\frac12}^{n}$\\\hline
&$x_{i-\frac12}$ &&$ x_{i+\frac12}$
\end{tabular} \end{center}

These are only 7 values (the box emphasizes that one of the values is a cell average, whereas the others are point values).

\subsubsection{Linear source term}

Consider first a linear source term, i.e. $s'' = 0$. Such source terms are relevant in practice (e.g. compressible Euler equations with gravity), and therefore it is worth dealing with them specifically as they allow for a simpler approach. For linear source it is possible to first find a quadrature for $q$ and to apply $s$ to the result. In order to find a quadrature formula for $q$, one needs to find a space-time polynomial $\mathscr P(t, x)$ of at least second degree which interpolates the available 7 data. Integrating this polynomial would yield a quadrature formula for $q$. Here we suggest to use
\begin{align}
 \mathscr P(t, x) = (a_0 + a_1 x + a_2 t + a_3 x^2 + a_4 x t + a_5 t^2) + a_6 x t^2 \label{eq:quadraturesource1d}
\end{align}
There is a unique set of coefficients $a_0, \ldots, a_6$ which makes polynomial \eqref{eq:quadraturesource1d} fulfill
\begin{align}
 \mathscr P(t^{n+1}, x_{i-\frac12}) &= q_{i-\frac12}^{n+1} & & & \mathscr P(t^{n+1}, x_{i+\frac12}) &= q_{i+\frac12}^{n+1} \label{eq:conditionquadrature1d1}\\
 \mathscr P(t^{n+\frac12}, x_{i-\frac12}) &= q_{i-\frac12}^{n+\frac12} & & & \mathscr P(t^{n+\frac12}, x_{i+\frac12}) &= q_{i+\frac12}^{n+\frac12} \\
 \mathscr P(t^{n}, x_{i-\frac12}) &= q_{i-\frac12}^n & \int\limits_{x_{i-\frac12}}^{x_{i+\frac12}} \!\!\! \dd x \, \mathscr P(t^n, x) &= q_i^n \Delta x & \mathscr P(t^n, x_{i+\frac12}) &= q_{i+\frac12}^{n}  \label{eq:conditionquadrature1d3}
\end{align}

Inserting this polynomial in \eqref{eq:numsource} and integrating it instead of the source yields the following quadrature formula:
\begin{align}
 \begin{split}
 \frac{1}{\Delta t}\int_0^{\Delta t} \dd t \, \frac{1}{\Delta x} \int_{-\frac{\Delta x}{2}}^{\frac{\Delta x}{2}} \dd x \, q(t^n + t, x_i + x) = \phantom{mmmmm}\\
 \bar q_i^n + \frac{1}{12} \left( - 5 (q_{i-\frac12}^n + q_{i+\frac12}^n) + q_{i-\frac12}^{n+1} + q_{i+\frac12}^{n+1} + 4 (q_{i-\frac12}^{n+\frac12} + q_{i+\frac12}^{n+\frac12}) \right ) 
 \end{split} \label{eq:quadrature1dlinearsource}
\end{align}
The weights can be depicted as

\setlength{\tabcolsep}{6.8pt}
\renewcommand{\arraystretch}{1.5}
\begin{center}\begin{tabular}{c|ccc}
$t^{n+1}$ &$\frac{1}{12}$  && $\frac{1}{12}$\\
 $t^{n+\frac12}$ & $\frac{4}{12}$ & & $\frac{4}{12}$\\
$t^n$ &$-\frac{5}{12}$ &$\boxed{1}$& $-\frac{5}{12}$\\\hline
& $x_{i-\frac12}$ && $x_{i+\frac12}$
\end{tabular}\end{center}
Again, the box indicates that the corresponding weight refers to the cell average, whereas the others multiply point values.

The time levels $(n, n+\frac12, n+1)$ contribute with weights $\left(\frac{1}{6}, \frac{2}{3}, \frac{1}{6}\right)$, such that this quadrature formula is a modification of Simpson's rule in time. Note that it is not possible to use terms proportional to $x^3$, $x^2 t$ or $t^3$ instead of the term $xt^2$ in the polynomial ansatz, as then the system \eqref{eq:conditionquadrature1d1}--\eqref{eq:conditionquadrature1d3} does not admit a solution. In a sense this is therefore the only choice of a simple quadrature formula.

Quadrature formula \eqref{eq:quadrature1dlinearsource} can be used immediately in order to approximate \eqref{eq:numsource} for linear source terms.

\subsubsection{Nonlinear source term}

For nonlinear $s$, the average 
\begin{align}
 \int\limits_{x_{i-\frac12}}^{x_{i+\frac12}} \!\!\! \dd x \, s(q(t^n, x)) 
\end{align}
in general is different from
\begin{align}
 s\left( \int\limits_{x_{i-\frac12}}^{x_{i+\frac12}} \!\!\! \dd x \, q(t^n, x) \right )
\end{align}

Point values, however, do not present any difficulties: one can just evaluate $s$ on them. Therefore we suggest to consider a reconstruction $q_{\text{recon},i}(x)$ that interpolates $q_{i-\frac12}^{n}$ and $q_{i+\frac12}^n$ and whose average agrees with $\bar q_i^n$. It is computed anyway in order to update the point values in time, see equation \eqref{eq:parabolicrecon}. This reconstruction can be easily evaluated at the midpoint of the cell. Then, instead of the cell averages, one works with a seventh point value $q_{\text{recon},i}(0) = \frac{1}{4} (6 \bar q_i^n - q_{i-\frac12}^n - q_{i+\frac12}^n)$. Of course, this is equivalent to replacing the average by a Simpson's rule in the quadrature, and thus the order of the quadrature is not reduced. Therefore when using only point values (the 6 pointwise degrees of freedom and one value at the cell midpoint) the weights of the quadrature formula read
\setlength{\tabcolsep}{6.8pt}
\renewcommand{\arraystretch}{1.5}
\begin{center}\begin{tabular}{c|ccc}
$t^{n+1}$ &$\frac{1}{12}$  && $\frac{1}{12}$\\
 $t^{n+\frac12}$ & $\frac{4}{12}$ & & $\frac{4}{12}$\\
$t^n$ &$-\frac{3}{12}$ &$ \frac8{12}$& $-\frac{3}{12}$\\\hline
& $x_{i-\frac12}$ && $x_{i+\frac12}$
\end{tabular}\end{center}

Equation \eqref{eq:numsource} then is replaced by the quadrature
{\footnotesize
\begin{align}
 \hat s_i = \frac{s(q^{n+1}_{i-\frac12}) + s(q^{n+1}_{i+\frac12}) + 4 \Big(s(q^{n+1}_{i-\frac12}) + s(q^{n+1}_{i+\frac12})\Big) - 3\Big(s(q^{n+1}_{i-\frac12}) + s(q^{n+1}_{i+\frac12})\Big) + 8 q_{\text{recon},i}(0) }{12} \label{eq:quadrature1dnonlinearsource}
\end{align}
}
This quadrature can now be used for nonlinear $s$. As \eqref{eq:quadrature1dnonlinearsource} uses a Simpson quadrature instead of the average, upon usage of a linear source $s$, it reduces to the expression \eqref{eq:quadrature1dlinearsource} because of the quadratic reconstruction.

If the source term vanishes, the scheme becomes conservative in the sense that averages are updated using numerical fluxes.

\subsection{Two spatial dimensions}

\subsubsection{Linear source term}

Similarly, consider the setup of the Active Flux method on two-dimensional Cartesian grids as described in \ref{ssec:dofs}. The available degrees of freedom are 
\begin{align*}
 3 \times 4 \text{ nodes: }&q_{i\pm\frac12,j\pm\frac12}^n, q_{i\pm\frac12,j\pm\frac12}^{n+\frac12}, q_{i\pm\frac12,j\pm\frac12}^{n+1}\\
 3 \times 2 \text{ vertical edges: }&q_{i\pm\frac12,j}^n, q_{i\pm\frac12,j}^{n+\frac12}, q_{i\pm\frac12,j}^{n+1}\\
 3 \times 2 \text{ horizontal edges: }&q_{i,j\pm\frac12}^n, q_{i,j\pm\frac12}^{n+\frac12}, q_{i,j\pm\frac12}^{n+1}\\
 1 \text{ average: }&\bar q_{ij}^n
\end{align*}

The ansatz for a space-time polynomial is
\begin{align}
 \mathscr P(t, x, y) = \left( \sum_{\zeta + \eta + \theta \leq 4} a_{\zeta\eta\theta} \cdot x^\zeta y^\eta t^\theta \right ) + a_{212} x^2 y t^2 + a_{122} x y^2 t^2
\end{align}

\begin{figure}
 \centering
 \includegraphics[width=0.5\textwidth]{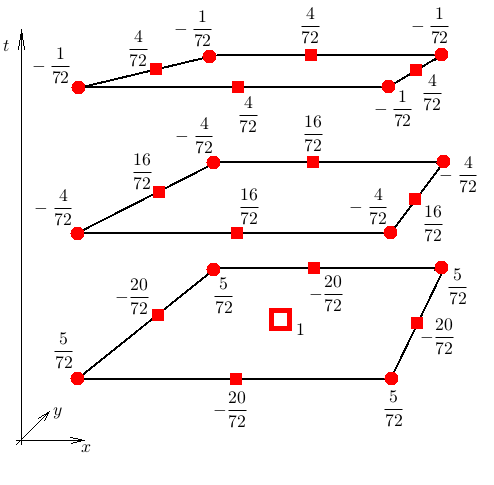}
 \caption{Illustration of the weights of the space time quadrature formula \eqref{eq:quadrature2dlinearsource}.}
 \label{fig:quadrature2d}
\end{figure}

It admits a unique solution to the interpolation problem given the available degrees of freedom and yields the following quadrature formula (see also figure \ref{fig:quadrature2d}):
\begin{align}
 \begin{split}
 \frac{1}{\Delta x}&\int_{-\frac{\Delta x}{2}}^{\frac{\Delta x}{2}} \dd x \, \frac{1}{\Delta y}\int_{-\frac{\Delta y}{2}}^{\frac{\Delta y}{2}} \dd y \, \frac{1}{\Delta t} \int_0^{\Delta t} \dd t \, q(t, x, y) = \bar q_{ij}^n  \\ 
 -&\frac{20}{72} \left( q_{\text E}^n
 + q_{\text N}^n
 + q_{\text S}^n
 + q_{\text W}^n \right )   + \frac{5}{72} \left( q_{\text NE}^n
 + q_{\text NW}^n
 + q_{\text SE}^n
 + q_{\text SW}^n \right ) \\
 +&\frac{16}{72} \left ( q_{\text E}^{n+\frac12} + q_{\text N}^{n+\frac12} + q_{\text S}^{n+\frac12} + q_{\text W}^{n+\frac12} \right )
 -\frac{4}{72} \left (q_{\text NE}^{n+\frac12} + q_{\text NW}^{n+\frac12} + q_{\text SE}^{n+\frac12} + q_{\text SW}^{n+\frac12} \right )\\
 +&\frac{4}{72} \left( q_{\text E}^{n+1} + q_{\text N}^{n+1} + q_{\text S}^{n+1} + q_{\text W}^{n+1} \right )
  - \frac{1}{72} \left( q_{\text NE}^{n+1} + q_{\text NW}^{n+1} + q_{\text SE}^{n+1} + q_{\text SW}^{n+1} \right ) 
  \end{split} \label{eq:quadrature2dlinearsource}
\end{align}
The time levels $(n, n+\frac12, n+1)$ contribute again with weights $(\frac16, \frac23, \frac16)$, and the edges always contribute $-4$ times the nodes.

\subsubsection{Nonlinear source term}

Again, for nonlinear source instead of the average it is necessary to use the evaluation of the reconstruction at the cell midpoint. This amounts to an approximation of the average by a two-dimensional Simpson rule. Then the source term is approximated as follows:
\begin{align}
 \begin{split}
 \frac{1}{\Delta x}\int_{-\frac{\Delta x}{2}}^{\frac{\Delta x}{2}} \dd x \, &\frac{1}{\Delta y}\int_{-\frac{\Delta y}{2}}^{\frac{\Delta y}{2}} \dd y \, \frac{1}{\Delta t} \int_0^{\Delta t} \dd t \, s(q(t, x, y)) = \frac{32}{72} s(q_{\text{recon},ij}(0, 0))  \\ 
 &-\frac{12}{72} \left( s(q_{\text E}^n)
 + s(q_{\text N}^n)
 + s(q_{\text S}^n)
 + s(q_{\text W}^n) \right )   \\&+ \frac{7}{72} \left( s(q_{\text NE}^n)
 + s(q_{\text NW}^n)
 + s(q_{\text SE}^n)
 + s(q_{\text SW}^n) \right ) \\
 &+\frac{16}{72} \left ( s(q_{\text E}^{n+\frac12}) + s(q_{\text N}^{n+\frac12}) + s(q_{\text S}^{n+\frac12}) + s(q_{\text W}^{n+\frac12}) \right )
 \\&-\frac{4}{72} \left (s(q_{\text NE}^{n+\frac12}) + s(q_{\text NW}^{n+\frac12}) + s(q_{\text SE}^{n+\frac12}) + s(q_{\text SW}^{n+\frac12}) \right )\\
 &+\frac{4}{72} \left( s(q_{\text E}^{n+1}) + s(q_{\text N}^{n+1} )+ s(q_{\text S}^{n+1}) + s(q_{\text W}^{n+1}) \right )
  \\&- \frac{1}{72} \left( s(q_{\text NE}^{n+1}) + s(q_{\text NW}^{n+1}) + s(q_{\text SE}^{n+1}) + s(q_{\text SW}^{n+1}) \right ) 
  \end{split} \label{eq:quadrature2dnonlinearsource}
\end{align}

In case that the data only depend on one of the spatial variables, the two-dimensional quadratures \eqref{eq:quadrature2dlinearsource} and \eqref{eq:quadrature2dnonlinearsource} do \emph{not} exactly reduce to the one dimensional quadratures \eqref{eq:quadrature1dlinearsource} and \eqref{eq:quadrature1dnonlinearsource}. This is because (cf. Figure \ref{fig:quadrature2d}) the point values on edge midpoints $\left(0, \pm\frac{\Delta y}{2}\right)$ do not disappear even if the data depend only on $x$, and therefore the available degrees of freedom remain different from the one-dimensional case.

\section{Well-balanced property for acoustics with gravity} \label{sec:wellbal}

\subsection{Exact evolution operator}

As described in \ref{ssec:evoacgrav}, a closed-form exact evolution operator for acoustics with gravity is very difficult to obtain. Nevertheless, it is still possible to show that a scheme endowed with such an operator would be well-balanced / stationarity preserving; i.e. that there exists a discretization of the stationary states of the PDE which remains exactly stationary. This proof does not require the evolution operator to be known explicitly, but only relies on the fact that the evolution operator is exact. Besides its fundamental importance, this result is used in section \ref{ssec:wbfix} to analyze the situation for the approximate evolution operator and to achieve the well-balanced property for it.

The numerical stationary states are best studied upon the (discrete) Fourier transform. Define $t_x := \exp(\ii k_x \Delta x)$, $t_y := \exp(\ii k_y \Delta y)$. Here $\ii$ is the imaginary unit and $\vec k = (k_x, k_y) \in \mathbb R^2$ is the wave vector characterizing the spatial frequency of the Fourier mode. Applying the Fourier transform introduces one mode $\bar q$ for the averages and one mode $q$ for the point values; this implies writing $q_{i} := \bar q t_x^i t_y^j$, $q_{i+\frac12} := q t_x^i t_y^j$.

\begin{theorem}[Stationarity preservation with exact evolution]
 If the discrete data fulfill
\begin{align}
  \bar \rho_i &= \frac{\rho_{i+\frac12} + \rho_{i-\frac12}}{2} \label{eq:acgravstatfinddiff0} \\
  \frac{ p_{i+ \frac12} - p_{i-\frac12} }{\Delta x} &=  g \frac{ \rho_{i-\frac12} + \rho_{i+\frac12} }{2 } \label{eq:acgravstatfinddiff1}\\
  \frac{\bar p_{i+1} - \bar p_{i} }{\Delta x} &=   g  \frac{ \rho_{i+\frac32} +4 \rho_{i+\frac12} +\rho_{i-\frac12} }{6 }  \label{eq:acgravstatfinddiff2}
\end{align}
 and the exact evolution operator for \eqref{eq:acgrav1}--\eqref{eq:acgrav3} is used, then the numerical solution remains stationary.
\end{theorem}
\begin{proof}
 The proof consists of two parts.
 \begin{enumerate}[i)]
  \item Consider first the evolution of the point values. When the exact evolution operator is used to update the point values, they remain stationary if the reconstruction fulfills 
\begin{align}
 v_\text{recon}(x) &= \mathrm{const} & \del_x p_\text{recon}(x) &= \rho_\text{recon}(x) g \label{eq:reconstationarity}
\end{align}

Upon the Fourier transform this becomes (w.l.o.g. $x_i = 0$)
\begin{align}
 - 3 \left(2 \bar p - p \left(1 + \frac{1}{t_x}\right)\right) \frac{2 x}{\Delta x^2} + p \left(1 - \frac{1}{t_x}\right) \frac{1}{\Delta x}  = \phantom{mmmmm} & \\\nonumber- 3 g \left(2 \bar \rho - \rho \left(1 + \frac{1}{t_x} \right)\right) \frac{x^2}{\Delta x^2} + g \rho \left(1 - \frac{1}{t_x}\right) \frac{x}{\Delta x} + g \frac{6\bar \rho - \rho \left(1 + \frac{1}{t_x}\right)}{4}
\end{align}
This shall be valid for all $x$:
\begin{align}
 2 \bar \rho - \rho (1 + 1/t_x) &= 0\\
 -2 \bar p t_x + p (t_x + 1)  &=  \frac{ \Delta x g \rho (t_x - 1) }{6}\\
 p (t_x - 1) &=  \Delta x g \frac{6\bar \rho t_x - \rho (t_x + 1)}{4}
\end{align}
These are three equations for four variables. In particular
\begin{align}
  \bar \rho &= \frac{\rho (1 + 1/t_x)}{2} \label{eq:acgravstatpointval1}\\
      p  &=  \Delta x g \rho \frac{  t_x + 1 }{2 (t_x-1)} \label{eq:acgravstatpointval2}\\
    \bar p  &=  \Delta x g \rho   \frac{  t_x^2 +4 t_x +1 }{6 t_x (t_x-1)} \label{eq:acgravstatpointval3}
\end{align}

These statements can be rewritten as finite difference formulae by inverting the Fourier transform to yield \eqref{eq:acgravstatfinddiff0}--\eqref{eq:acgravstatfinddiff2}.

\item Assume now \eqref{eq:acgravstatpointval1}--\eqref{eq:acgravstatpointval3} to be true. Simpson's rule in time for the flux average is trivial, and thus the update of the cell average amounts to
\begin{align}
 \frac{\bar v^{n+1} - \bar v^n}{\Delta t} + \frac{p (1 - 1/t_x)}{\Delta x} &= \frac{\bar v^{n+1} - \bar v^n}{\Delta t} +  g \rho \frac{  t_x + 1 }{2 t_x}\\
 &= \frac{\bar v^{n+1} - \bar v^n}{\Delta t} +  g \bar \rho 
\end{align}

The quadrature formula \eqref{eq:quadrature1dlinearsource} for the source reduces to $g \bar \rho $ if the point values are stationary, which implies $\bar v^{n+1} = \bar v^n$. This completes the proof.

 \end{enumerate}
\end{proof}

The equations \eqref{eq:acgravstatpointval1}--\eqref{eq:acgravstatpointval3} contain $\rho$ as a free variable. One can rewrite the system making $p$ the free variable:
\begin{align}
 \bar \rho &= \frac{p (t_x-1)}{t_x \Delta x g} &
 \rho &= \frac{2 p (t_x-1) }{\Delta x g (t_x+1)} &
 \bar p &= p\frac{t_x^2 + 4 t_x + 1}{3 t_x (t_x + 1)} \label{eq:acgravstatpointval4}
\end{align}
This form will be useful later.

Equations \eqref{eq:acgravstatfinddiff1}--\eqref{eq:acgravstatfinddiff2} are finite difference approximations of $\del_x p = \rho g$. {By construction, the discrete stationary states are those whose reconstruction fulfills \eqref{eq:reconstationarity} in every cell.} Equation \eqref{eq:acgravstatfinddiff0} implies that the reconstructed $\rho$ of the discrete stationary state is linear, which is clear: for quadratic reconstructions to fulfill \eqref{eq:reconstationarity}, $\rho_\text{recon}$ has to be linear in each cell. The slope of the linear function can vary from cell to cell.

\subsection{Approximate evolution operator} \label{ssec:wbfix}

The above section identifies conditions \eqref{eq:acgravstatfinddiff0}--\eqref{eq:acgravstatfinddiff2} on the discrete data for them to remain stationary upon usage of the \emph{exact} evolution operator. Unfortunately, such an operator is unavailable in practice. Having identified an approximate solution operator, which agrees with the exact solution up to terms $\mathcal O(t^3)$ in section \ref{ssec:evolinhyp}, here we study whether it keeps the same data \eqref{eq:acgravstatfinddiff0}--\eqref{eq:acgravstatfinddiff2} stationary as well.

\begin{theorem}
 If the discrete data fulfill \eqref{eq:acgravstatfinddiff0}--\eqref{eq:acgravstatfinddiff2} and the approximate evolution operator of theorem \ref{thm:rk2} for \eqref{eq:acgrav1}--\eqref{eq:acgrav3} is used, then both the pressure $p$ and the density $\rho$ remain stationary over one time step, but the velocity undergoes the time evolution
 \begin{align}
  v_{i+\frac12}(t) &= - \frac{\alpha g^2}{4}  \frac{ \rho_{i+\frac12} - \rho_{i-\frac12}  }{\Delta x}  t^3 \label{eq:vstarrewrite}
 \end{align}
\end{theorem}
\begin{proof}
Assume the initial data to fulfill \eqref{eq:acgravstatfinddiff0}--\eqref{eq:acgravstatfinddiff2}, or equivalently \eqref{eq:reconstationarity}. Using \eqref{eq:parabolicrecon} (and applying the discrete Fourier transform straight away) \eqref{eq:reconstationarity} implies
\begin{align}
 p_\text{recon}(x) &= \frac14 \left(6 \bar p - p \left(1 + \frac{1}{t_x}\right)\right) + \frac{x}{\Delta x} \left(1 - \frac{1}{t_x}\right ) p - 3\frac{x^2}{\Delta x^2} \left( 2 \bar p - p \left(1 + \frac{1}{t_x} \right )  \right )\\
 \rho_\text{recon}(x) &= \frac{1}{g \Delta x} \left( p \left( 1 - \frac{1}{t_x}  \right ) - 6 \frac{x}{\Delta x} \left( 2 \bar p - p \left(1 + \frac{1}{t_x} \right )   \right )  \right )\\
 v_\text{recon}(x) &= 0
\end{align}
and using \eqref{eq:acgravcharvar} therefore
\begin{align}
 Q_{1,0}(x) =\; & Q_{2,0}(x) = -\frac{p(1+  t_x)-6 {\bar p}  t_x  }{8  t_x }+\frac{p (t_x -1) x}{2 \Delta x   t_x }+\frac{3
(p(1+  t_x)-2 {\bar p}  t_x  ) x^2}{2 \Delta x ^2  t_x }\\
 Q_{3,0}(x) =\; & \frac{ p (-1+ t_x )}{\Delta x  g t_x} + \frac{p-6 {\bar p}  t_x +p  t_x }{4 c^2 t_x} \\&+\frac{\left(-\Delta x 
g p (t_x -1)+6 c^2 (p(1+  t_x)-2 {\bar p}  t_x  )\right) x}{c^2 \Delta x ^2 g  t_x }-\frac{3 (p(1+ t_x)-2 {\bar p}  t_x  ) x^2}{c^2 \Delta x ^2  t_x }\nonumber
\end{align}

Evaluating the Runge-Kutta algorithm of section \ref{ssec:evolinhyp} on these initial data (at $x = \frac{\Delta x}{2}$) yields
\begin{align}
 (\rho ,v^*, p)^\text T \qquad
\text{with}
  \qquad v^* = - \frac{\alpha g (t_x-1)^2}{2 \Delta x^2 t_x (t_x+1)} p t^3 \label{eq:vstar}
\end{align}
($\alpha$ is the parameter appearing in the RK2 method.)

Recall that $\rho$ and $p$ are the Fourier coefficients of the point values of the density and the pressure. Obviously $\rho$ and $p$ remain stationary, but the velocity does not. Using \eqref{eq:acgravstatpointval4} $v^*$ can be rewritten as
\begin{align}
 v^* &= - \frac{\alpha g^2}{4 \Delta x}  \left( 1 - \frac{1}{t_x}    \right )   \rho t^3 = - \frac{\alpha g^2}{4}  \frac{ \rho_{i+\frac12} - \rho_{i-\frac12}  }{\Delta x}  t^3 
\end{align}
having applied the inverse Fourier transform in the last step.
\end{proof}

Observe that the time evolution of the velocity is consistent with the accuracy of the algorithm ($\mathcal O(t^3)$).

\begin{corollary}[Stationarity preservation with approximate evolution]
If the algorithm of section \ref{ssec:evolinhyp} is modified by adding the term
\begin{align}
 \frac{\alpha g^2}{4}  \frac{ \rho_{i+\frac12} - \rho_{i-\frac12}  }{\Delta x}  t^3 \label{eq:velrk2modifcation}
\end{align}
to the velocity evolution, then 
\begin{enumerate}[i)]
 \item its accuracy is not changed
 \item it becomes stationarity preserving / well-balanced with the same discrete stationary states as the exact evolution operator.
\end{enumerate}
\end{corollary}

The two forms \eqref{eq:vstar} and \eqref{eq:vstarrewrite} of $v^*$ are equivalent, because the initial data have been chosen to be stationary, and thus additionally fulfill \eqref{eq:acgravstatpointval4}. The proposed modification is to \emph{always} add $-v^*$ to the velocity evolution, irrespective of whether the data fulfill \eqref{eq:acgravstatpointval4} or not. At this point the Fourier coefficients of $\rho$ and $p$ are independent and it matters whether the correction is used in the form \eqref{eq:vstar} or \eqref{eq:vstarrewrite}. Of course, also the inverse Fourier transform has to be applied to the expression first in order for the correction to attain the form of a finite difference formula. Compact finite difference formulae are in one-to-one-correspondence with Laurent polynomials in $t_x$. An expression such as $\frac{1}{t_x + 1} = 1 - t_x + t_x^2 \mp \ldots$ is an expression involving an unbounded stencil and cannot be implemented in usual codes. Therefore \eqref{eq:vstarrewrite} cannot be used as a correction because the correction would have a non-compact stencil (just as the equivalent expressions involving only $\bar \rho$ or $\bar p$). This is why the form \eqref{eq:vstar} which involves point values of $\rho$ is preferred. 

Being always present in the velocity evolution (and not only at stationary states), the modification \eqref{eq:velrk2modifcation} might in general affect the stability of the algorithm, but it has not been found to have any effect on the stability in practice.

\section{Numerical examples} \label{sec:numeric}

The numerical examples of this section serve to illustrate the performance of the new method. The equations discussed are linear advection with different source terms (in one and two spatial dimensions, as introduced in section \ref{ssec:evolinadv}) and linear acoustics with gravity (introduced in section \ref{ssec:evoacgrav}). In both cases it is demonstrated that the method achieves third order of accuracy in the experiments. For acoustics with gravity additionally the discrete stationary states are studied and shown to agree with the prediction of section \ref{sec:wellbal}.

\subsection{Linear advection}

\begin{figure}
 \centering
 \includegraphics[width=0.48\textwidth]{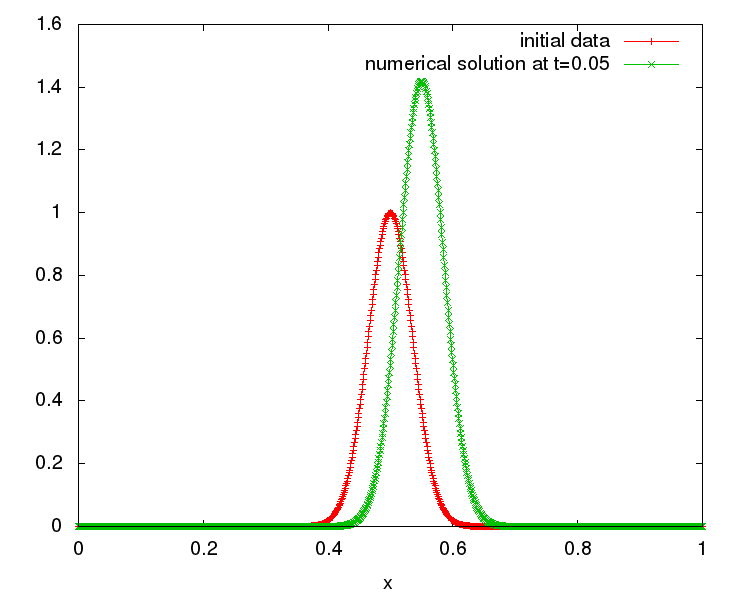} 
 \hfill
 \includegraphics[width=0.48\textwidth]{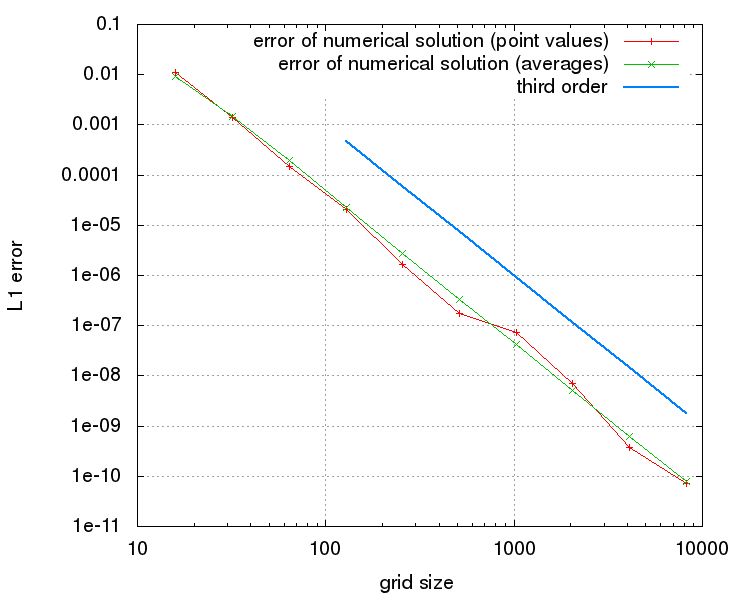} 
 \caption{Gaussian initial data for \eqref{eq:linadvlinsource} with $\vec U = \vec e_x$, $\kappa = 7$. Note that due to the source term, the Gaussian is advected and also changes shape. Exact evolution operator \eqref{eq:linadvlinsourceexact} and quadrature formula \eqref{eq:quadrature1dlinearsource} have been used with CFL = {0.9}. \emph{Left}: Initial data and solution at $t=0.05$ (cell averages) on a grid with 1000 cells. \emph{Right}: Error of the numerical solution as a function of the grid size shows third order convergence.}
 \label{fig:advectionexact}
\end{figure}

\begin{figure}
 \centering
 \includegraphics[width=0.48\textwidth]{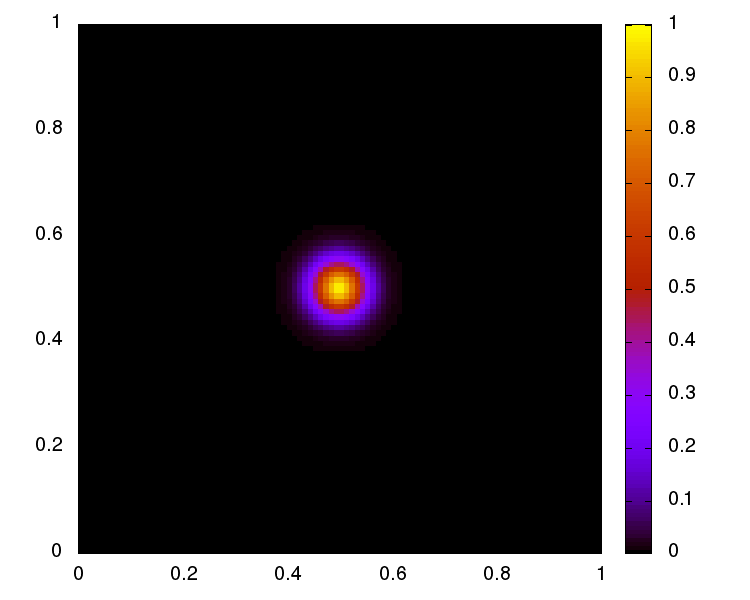} 
 \hfill
 \includegraphics[width=0.48\textwidth]{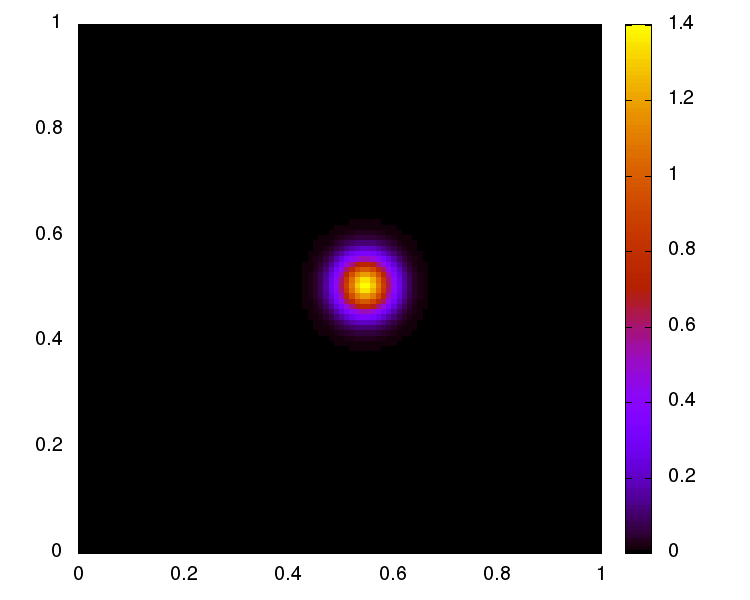} 
 \caption{Gaussian initial data for \eqref{eq:linadvlinsource} with $\vec U = (1, 0.1)$, $\kappa = 7$. Note that due to the source term, the Gaussian is advected and also changes shape. Exact evolution operator \eqref{eq:linadvlinsourceexact} and quadrature formula \eqref{eq:quadrature2dlinearsource} have been used with CFL = {0.9}. \emph{Left}: Initial setup. \emph{Right}: Numerical solution at $t=0.05$ on a $100 \times 100$ Cartesian grid.}
 \label{fig:advectionexactmultidsetup}
\end{figure}

\begin{figure}
 \centering
 \includegraphics[width=0.48\textwidth]{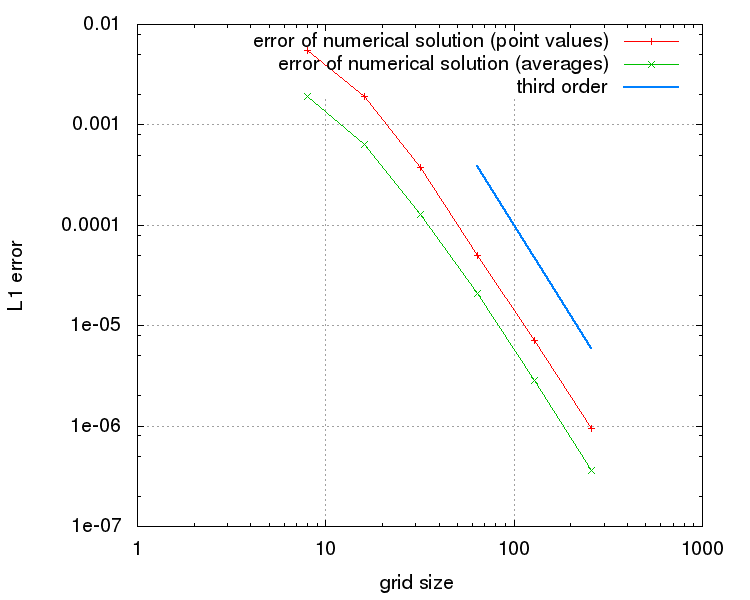} 
 \caption{Convergence study for the setup shown in Figure \ref{fig:advectionexactmultidsetup}. One observes third order accuracy.}
 \label{fig:advectionexactmultid}
\end{figure}

Consider first 
\begin{align}
 \del_t q + \vec U \cdot \nabla  q &= \kappa q \label{eq:linadvlinsource}
\end{align}
with the exact solution given by \eqref{eq:linadvlinsourceexact}. In Figures \ref{fig:advectionexact}--\ref{fig:advectionexactmultid} the exact solution operator is used for the evolution of the point values and third order convergence is observed. This shows that the quadrature formulae \eqref{eq:quadrature1dlinearsource} and \eqref{eq:quadrature2dlinearsource} used to evolve the cell averages indeed yield a third order scheme. Figure \ref{fig:advectionexact} shows the setup for a one-dimensional situation together with a convergence study, Figure \ref{fig:advectionexactmultidsetup} shows the setup in two spatial dimensions and Figure \ref{fig:advectionexactmultid} shows the corresponding convergence study.


\begin{figure}
 \centering
 \includegraphics[width=0.48\textwidth]{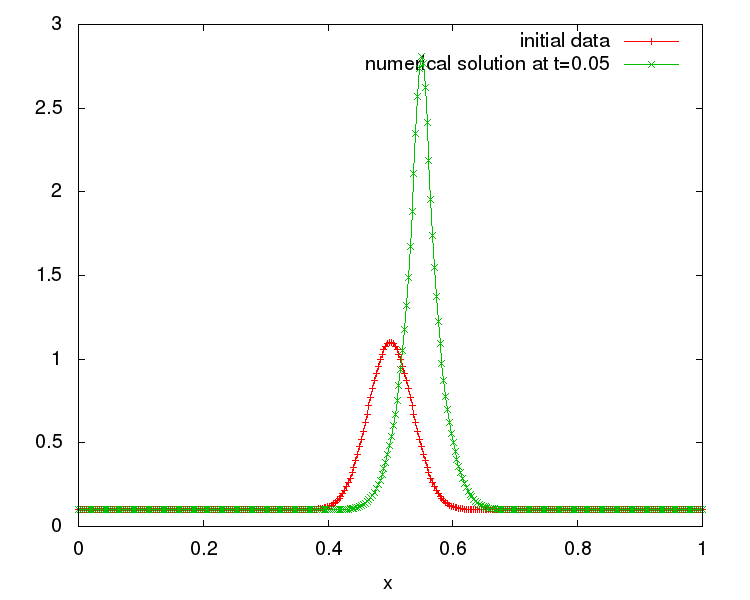} 
 \hfill
 \includegraphics[width=0.48\textwidth]{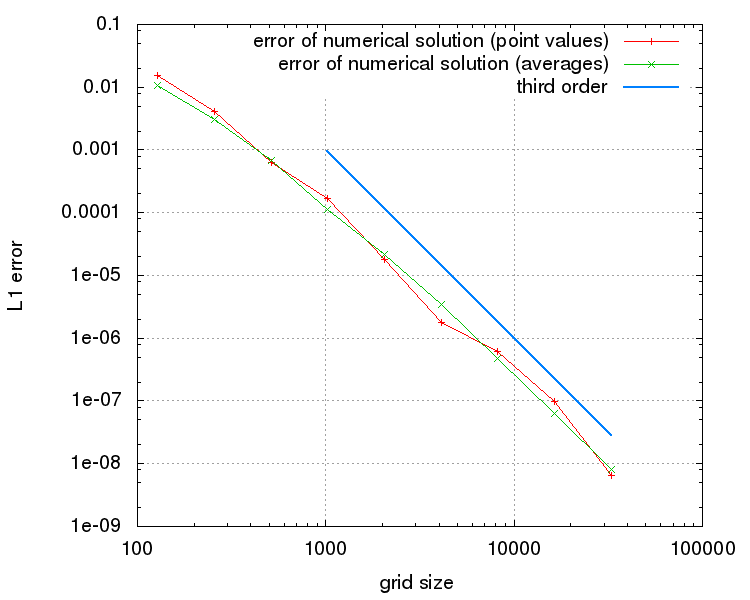} 
 \caption{Gaussian initial data for \eqref{eq:linadvnonsource} with $s(q) = \kappa q^B$ and $\vec U = \vec e_x$, $\kappa = 7$, $B = 3$. Runge-Kutta approximate evolution operator from Corollary \ref{cor:scalarrk2} (with $\alpha = \frac12$) and quadrature formula \eqref{eq:quadrature1dnonlinearsource} have been used with CFL = {0.9}. The solution has been computed on a grid covering $[-1:2]$, but the error is only computed inside $[0,1]$ to exclude any boundary influence. \emph{Left}: Initial setup and solution at $t=0.05$ (cell averages) on a grid with 1000 cells. \emph{Right}: Error of the numerical solution as a function of the grid size shows third order convergence. The exact solution is given by \eqref{eq:linadvnonlinsourceexact}.}
 \label{fig:advectionnonlin}
\end{figure}

Consider now
\begin{align}
 \del_t q + \vec U \cdot \nabla  q &= \kappa q^B \qquad B \neq 1 \label{eq:linadvnonsource}
\end{align}
with the exact solution \eqref{eq:linadvnonlinsourceexact} and $\kappa = 7$, $B = 3$. Figure \ref{fig:advectionnonlin} (left) shows the initial data and the numerical solution, and Figure \ref{fig:advectionnonlin} (right) shows a convergence study for the approximate evolution operator from Corollary \eqref{cor:scalarrk2}. One observes third order accuracy, as expected.

\subsection{Acoustics with gravity}

Consider now the equations of linear acoustics with a gravity source term \eqref{eq:acgrav1}--\eqref{eq:acgrav3}. The exact solution operator is only partly available in closed form, and therefore the approximate Runge-Kutta evolution operator of section \ref{ssec:evolinhyp} is used in combination with the well-balancing fix \eqref{eq:velrk2modifcation}. The parameter $\alpha$ in the Runge-Kutta method is chosen to $\alpha = \frac12$ and CFL = {0.9} everywhere.

\begin{figure}
 \centering
 \includegraphics[width=0.48\textwidth]{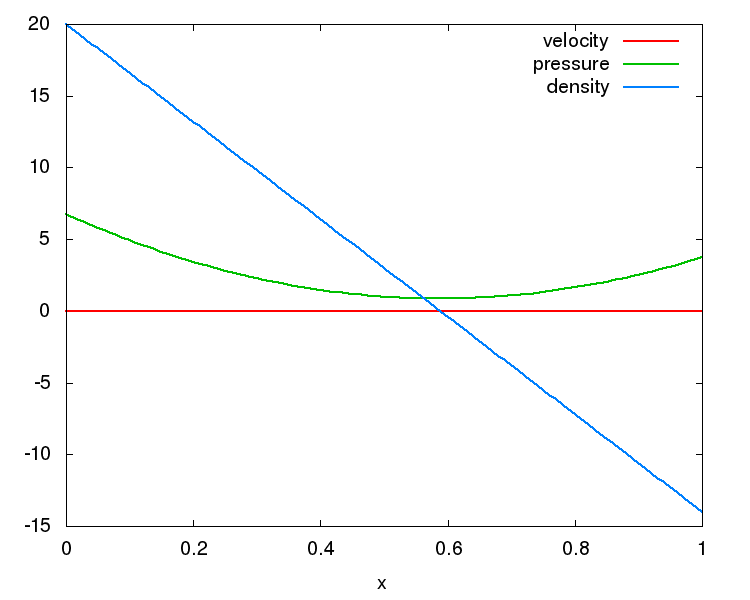} 
 \hfill
 \includegraphics[width=0.48\textwidth]{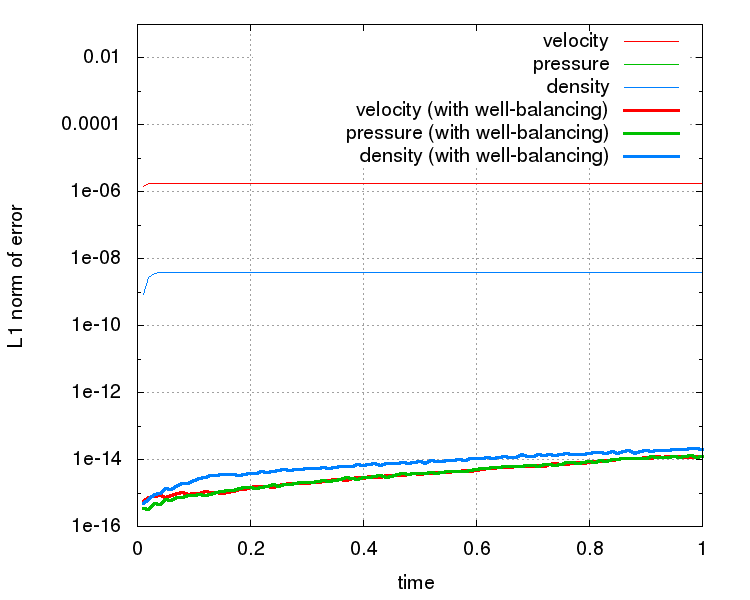} 
 \caption{Setup of a stationary parabola \eqref{eq:statparabolasetup} for \eqref{eq:acgrav1}--\eqref{eq:acgrav3}, solved using the Runge-Kutta approximate evolution operator of section \ref{ssec:evolinhyp} with and without well-balancing \eqref{eq:velrk2modifcation}. Here $g = -1$, and the setup is solved on a grid covering $[-1.5,2.5]$, but the error is only measured inside $[0,1]$ ($\Delta x = 10^{-2}$) to exclude the influence of the boundaries. \emph{Left}: Setup. \emph{Right}: Error of numerical solution (point values) as a function of time. Thin lines: without the well-balancing \eqref{eq:velrk2modifcation}. Thick lines: including the well-balancing \eqref{eq:velrk2modifcation}. In the latter case one only observes an evolution due to machine error.}
 \label{fig:acgravstatparabola}
\end{figure}

Figure \ref{fig:acgravstatparabola} shows a stationary setup given by
\begin{align}
 p &= A_1 x^2 + A_2x + A_3 & \rho &= 2A_1x/g + A_2/g & v &= 0 \label{eq:statparabolasetup}
\end{align}
with $A_1 = 17, A_2 = -3, A_3 = 1$. This parabola is exactly recovered by the reconstruction, and thus remains stationary up to machine precision. This experiment shows that the well-balancing fix works as it should.

\begin{figure}
 \centering
 \includegraphics[width=0.48\textwidth]{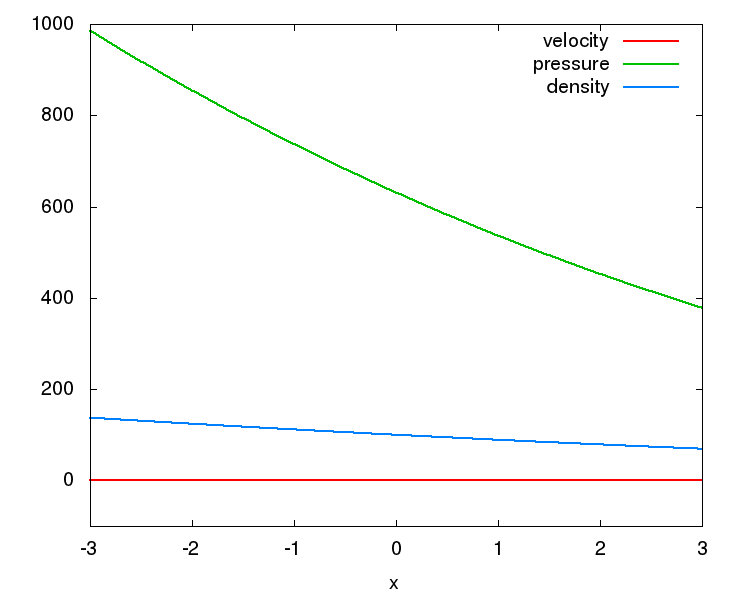} 
 \hfill
 \includegraphics[width=0.48\textwidth]{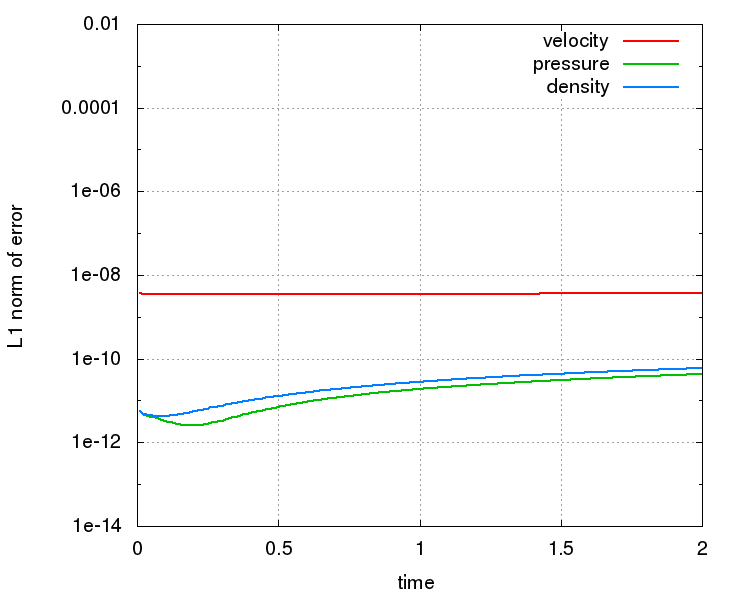} 
 \caption{Stationary setup \eqref{eq:statatmosetup} for \eqref{eq:acgrav1}--\eqref{eq:acgrav3}, solved using the Runge-Kutta approximate evolution operator of section \ref{ssec:evolinhyp} with well-balancing \eqref{eq:velrk2modifcation}. Here $g = -1$, and the setup is solved on a grid covering $[-5.5,5.5]$, but the error is only measured inside $[-3,3]$ ($\Delta x = 1/300$) to exclude the influence of the boundaries. \emph{Left}: Setup (cell averages). \emph{Right}: Error of numerical solution (point values) as a function of time. One observes a transition towards a numerical stationary state which then persists forever.}
 \label{fig:statatmo}
\end{figure}

Consider next (Figure \ref{fig:statatmo}) the stationary setup fulfilling $p = K \rho^\gamma$, i.e.
\begin{align}
 \rho = \left ( \frac{g(\gamma-1) }{K\gamma}   x + \rho_0^{\gamma-1} \right )^{\frac1{\gamma-1}}  \label{eq:statatmosetup}
\end{align}
with $K=1, \gamma = 1.4$, $\rho_0 = 100$. This is reminiscent of an isentropic atmosphere in the context of the Euler equations. This setup is not recovered exactly by the reconstruction, but one observes a numerical evolution towards a discrete stationary state which then persists forever.

\begin{figure}
 \centering
 \includegraphics[width=0.48\textwidth]{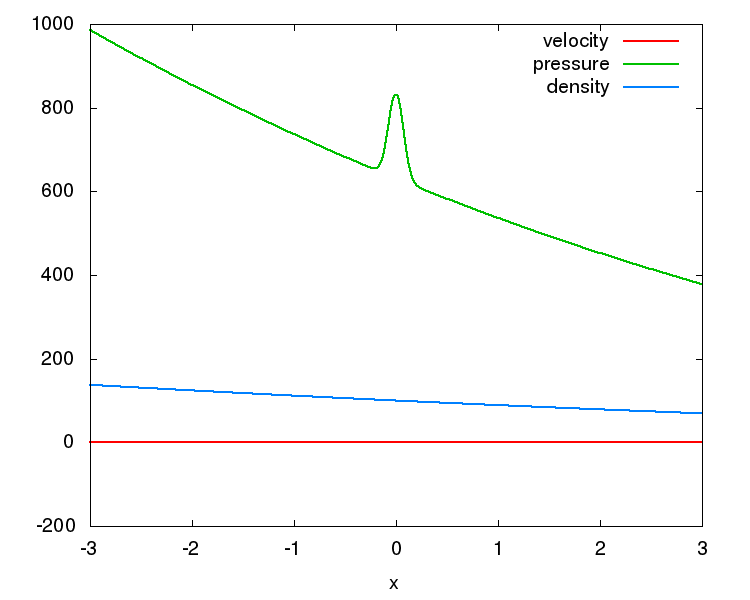} 
 \hfill
 \includegraphics[width=0.48\textwidth]{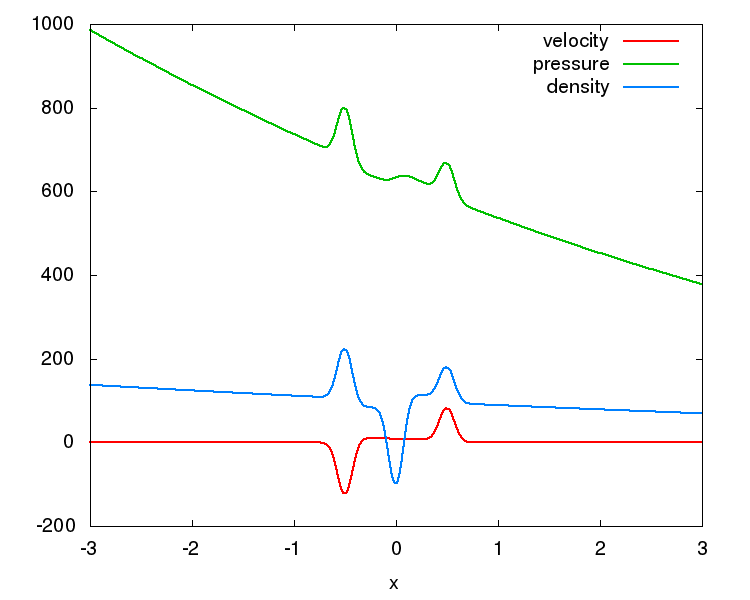} 
 \caption{Setup \eqref{eq:statatmosetup} endowed with the pressure perturbation \eqref{eq:statatmopert} solved using the Runge-Kutta approximate evolution operator of section \ref{ssec:evolinhyp} with well-balancing \eqref{eq:velrk2modifcation}. \emph{Left}: Initial data (cell averages). \emph{Right:} Numerical solution (cell averages) at $t=0.5$ on a grid covering $[-5.5, 5.5]$, but only the subinterval $[-3,3]$ is considered in order to exclude the influence of the boundaries. $\Delta x = 0.01$, CFL = {0.9}.}
 \label{fig:statatmopert}
\end{figure}

\begin{figure}
 \centering
 \includegraphics[width=0.48\textwidth]{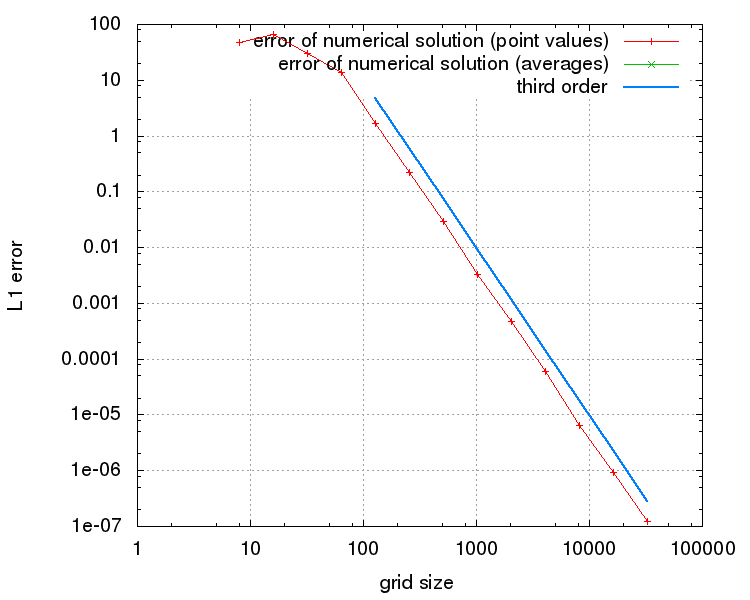} 
 \caption{Setup of Figure \ref{fig:statatmopert}. The error of the numerical solution is measured on the point values. One observes third order accuracy.}
 \label{fig:statatmopertconv}
\end{figure}

Next, a perturbation
\begin{align}
 200 \exp(- 100 x^2) \label{eq:statatmopert}
\end{align}
in the pressure is added onto the setup \eqref{eq:statatmosetup}. In order to study the accuracy of the scheme on this setup, it is solved on a grid of $131072 = 2^{18}$ cells and the solution is used as reference. Again, $g=-1, K = 1, \gamma = 1.4$. Figure \ref{fig:statatmopert} shows the setup and the numerical solution at $t=0.5$, and Figure \ref{fig:statatmopertconv} shows a convergence study which displays third order convergence.


\begin{figure}
 \centering
 \includegraphics[width=0.48\textwidth]{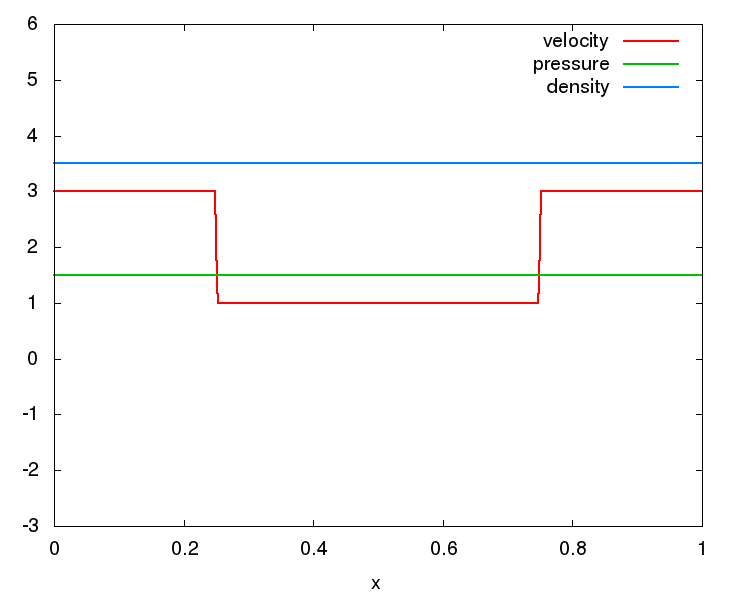} 
 \hfill
 \includegraphics[width=0.48\textwidth]{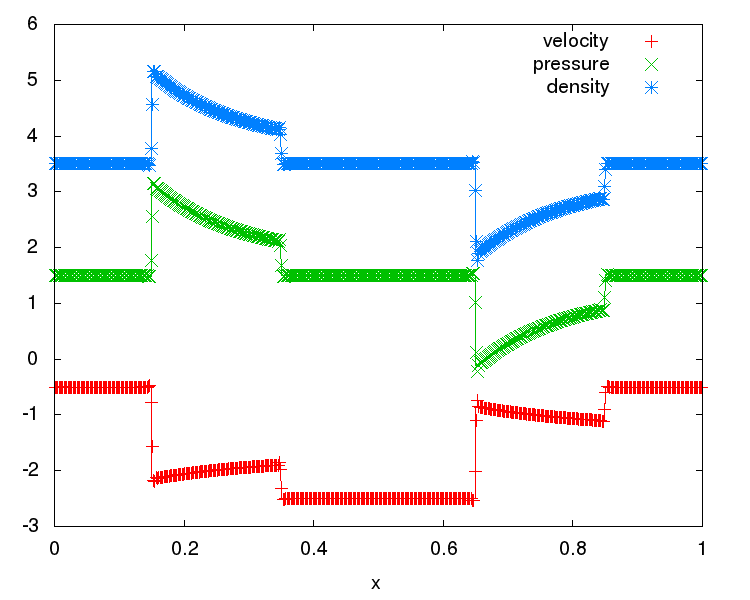} 
 \caption{Riemann problem setup \eqref{eq:acgravriemannproblem} solved using the Runge-Kutta approximate evolution operator of section \ref{ssec:evolinhyp} with well-balancing \eqref{eq:velrk2modifcation}. Here, $g=-10$. \emph{Left}: Initial data. \emph{Right:} Numerical solution (dots) and exact solution (solid line) at $t=0.1$. $\Delta x = 0.01$, CFL = {0.9}. Averages of the numerical solution are shown are shown.}
 \label{fig:ac-rp}
\end{figure}

Consider finally a Riemann problem:
\begin{align}
\rho &= 3.5 &  p &= 1.5 & v &= \begin{cases} 1 & 0.25 \leq x \leq 0.75 \\ 3 & \text{else}  \end{cases} \label{eq:acgravriemannproblem}
\end{align}
This Riemann problem can be solved exactly using the formula \eqref{eq:exactsolacgravV1}--\eqref{eq:exactsolacgravV2}. Note that if all quantities are constant in space, then they solve
\begin{align}
 \del_t \rho &= 0 &
 \del_t p &= 0 &
 \del_t v &= \rho g
\end{align}
which means that $\rho$ and $p$ remain stationary, but that $v = v(t=0) + \rho g t$. The solution to the initial data \eqref{eq:acgravriemannproblem} therefore can be obtained by adding the time evolution of $(0,v_0(x),0)^\text{T}$ (via numerical quadrature of \eqref{eq:exactsolacgravV1}--\eqref{eq:exactsolacgravV2}) and the time evolution of $(\rho,0,p)^\text T$ which is just $(\rho,\rho gt,p)^\text{T}$. Figure \ref{fig:ac-rp} shows the numerical and the exact solution.

\section{Conclusions and outlook}

Active flux is a novel kind of numerical method for hyperbolic problems, extending the finite volume method. Instead of computing the intercell flux via a Riemann problem it relies on a continuous reconstruction and on accurately evolved point values along the cell boundary. They then immediately serve as quadrature values for the computation of the intercell flux. The extension of Active Flux to time dependent balance laws presented in this paper requires a modification in both these aspects: the evolution of the point values and the average update need to account for the source term. Here, an approximate evolution operator is suggested for the point value update; this is done for linear systems with possibly nonlinear source terms in one spatial dimension, and linear scalar equations with source terms in multiple spatial dimensions. A suitable quadrature is suggested in order to approximate the contribution of the source term to the cell average. This quadrature can be applied to any system of (nonlinear) balance laws. 

We aim at combining the strategy presented in this paper with an approximate evolution operator for a nonlinear homogeneous problem (such as those suggested in \cite{barsukow19activeflux}) in future. Multi-dimensional systems of hyperbolic conservation laws are very different from their one-dimensional counterparts because in general characteristics are unavailable and need to be conceptually replaced by characteristic cones. Examples of evolution operators that make use of such cones can be found in \cite{eymann13,fan15,fan17,barsukow18activeflux}. Combining these with an approximate evolution of the source term shall pave the way towards the extension of Active Flux to nonlinear multi-dimensional balance laws and the derivation of accurate structure preserving (in particular well-balanced) methods for them.

\section*{Acknowledgement}

We thank Philip L. Roe for valuable comments and advice. WB was supported by the German Academic Exchange Service (DAAD) with funds from the German Federal Ministry of Education and Research (BMBF) and the European Union (FP7-PEOPLE-2013-COFUND -- grant agreement no. 605728) as well as by the Deutsche Forschungsgemeinschaft (DFG) through project 429491391 (BA 6878/1-1).

\newcommand{\etalchar}[1]{$^{#1}$}

\appendix

\section{Exact solution of linear acoustics with gravity} \label{sec:excatsolutionac}

System \eqref{eq:acgrav1}--\eqref{eq:acgrav3} can in principle be immediately solved exactly via Fourier transform by inserting the ansatz 
\begin{align}
 \veccc{\rho}{v}{p} = \veccc{\hat\rho}{\hat v}{\hat p} \exp(\ii k \cdot x - \ii \omega t)
\end{align}
into \eqref{eq:acgrav1}--\eqref{eq:acgrav3}:
\begin{align}
   \omega \veccc{\hat\rho}{\hat v}{\hat p} &= \left( \begin{array}{ccc} 0 &  k & 0 \\ \ii g & 0 & k \\ 0 &  c^2 k & 0 \end{array} \right ) \veccc{\hat\rho}{\hat v}{\hat p} 
\end{align}
Therefore $\omega = 0$, or $\omega = \pm\sqrt{c^2 k^2 + \ii g k}$. The complex eigenvalue can be removed upon transforming
\begin{align}
 \rho &= \tilde \rho \ee^{\mu x} & v &= \tilde v \ee^{\mu x} & p &= \tilde p \ee^{\mu x} 
\end{align}
with
\begin{align}
 \mu := \frac{g}{2 c^2}
\end{align}

System \eqref{eq:acgrav1}--\eqref{eq:acgrav3} then reads
\begin{align}
 \del_t \tilde \rho + \del_x \tilde v &= -\mu \tilde v \label{eq:acgravtransform1}\\
 \del_t \tilde v + \del_x \tilde p &= \tilde \rho g - \mu \tilde p \\
 \del_t \tilde p + c^2 \del_x \tilde v &= - c^2 \mu \tilde v \label{eq:acgravtransform3}
\end{align}
Now, a solution of \eqref{eq:acgravtransform1}--\eqref{eq:acgravtransform3} shall be found. For better readability, drop the tilde. Upon the Fourier transform \eqref{eq:acgravtransform1}--\eqref{eq:acgravtransform3} becomes
\begin{align}
   \omega \veccc{\hat\rho}{\hat v}{\hat p} &= \mathcal E \veccc{\hat\rho}{\hat v}{\hat p} & \mathcal E &= \left( \begin{array}{ccc} \phantom{m}0\phantom{ml} &  k-\ii \mu & 0 \\ \ii g\phantom{il} & 0 & k-\ii \mu \\ \phantom{m}0\phantom{ml} &  c^2 k-\ii c^2 \mu & 0 \end{array} \right )
\end{align}

The eigenvalues of $\mathcal E$ are now real: $\omega_1 = 0$, $\omega_{2,3} = \pm c \sqrt{ k^2 + \mu^2} $. Although this transformation brings the endeavor of finding the exact solution to \eqref{eq:acgrav1}--\eqref{eq:acgrav3} into the realm of the possible, technical difficulties prevent one from actually computing all Green's functions in closed form. 

Assume therefore that the only non-vanishing initial data are in the velocity. Then the Fourier mode at initial time reads
\begin{align}
 (0,\hat v,0)^\text{T} \exp(\ii k x) \label{eq:initialvfouriermode}
\end{align}
and at a later time it becomes
\begin{align}
  \sum_{m=1}^3 v_m \exp(\ii k x - \ii \omega_m t)
\end{align}
where the decomposition of $(0,\hat v,0)^\text{T}$ in the eigenbasis of $\mathcal E$ is used, i.e.
\begin{align}
 (0,\hat v,0)^\text{T} &= \sum_{m=1}^3 v_m & \mathcal E v_m &= \omega_m v_m
\end{align}
Such a basis is given e.g. by
\begin{align}
 e_1 &= \veccc{\mu + \ii  k }{0}{ g} &
 e_{2,3} &= \veccc{\mu + \ii k}{\pm  \ii c\sqrt{k^2 + \mu^2} }{c^2 (\mu +  \ii k)} 
 \end{align}
Collecting the terms yields the time evolution of the Fourier mode \eqref{eq:initialvfouriermode}:
\begin{align}
 &\hat v \exp(\ii k x) \veccc{\displaystyle - \frac{(\mu + \ii k) \sin \left( ct\sqrt{k^2 + \mu^2}  \right ) }{c \sqrt{k^2 + \mu^2} }  }{\displaystyle  \cos \left( ct \sqrt{k^2 + \mu^2}   \right )   }{\displaystyle   - \frac{c^2 (\mu + \ii k) \sin \left( ct\sqrt{k^2 + \mu^2}  \right ) }{ c\sqrt{k^2 + \mu^2} }    }\\
&= \hat v \veccc{\displaystyle - (\mu + \del_x)   }{\displaystyle  \del_t \phantom{\frac12}   }{\displaystyle   - c^2(\mu + \del_x  )    }      \exp(\ii k x) \frac{\sin \left(c t\sqrt{k^2 + \mu^2}  \right )}{c\sqrt{k^2 + \mu^2}}
\end{align}

Green's function is obtained by inserting the Fourier transform of a Dirac $\delta_{x'}$ at $x'$, i.e. taking $\hat v = \frac{\exp(-\ii k x')}{\sqrt{2\pi}}$ and performing the inverse Fourier transform with the help of formula 1.7\,(30) in \cite{bateman54}. This yields, wherever defined,
\begin{align}
 \veccc{\displaystyle G_\rho(t, x; x')}{\displaystyle G_v(t, x; x')}{\displaystyle  G_p(t, x; x')} &= 
 \veccc{\displaystyle 
 - (\mu + \del_x  )  
 }{\displaystyle  
 \del_t \phantom{\frac12}   }{\displaystyle   
 - c^2 (\mu + \del_x  )    }  
 \frac1{2c} J_0\left(\mu \sqrt{(ct)^2 - (x-x')^2}\right)
 \\\nonumber&+ \veccc{\displaystyle - \frac{\delta_{x+ct} - \delta_{x-ct}}{2c} }{\displaystyle \frac{\delta_{x+ct} + \delta_{x-ct}}{2}  }{\displaystyle c \left(\delta_{x+ct} - \delta_{x-ct}\right)}
\end{align}
where $J_0$ is the 0-th order Bessel function of the first kind, and $J_0' = -J_1$. Then the solution is obtained by performing a convolution with the initial data. Reinstalling the tilde one has
\begin{align}
 \tilde v(t, x) &= \int \dd x'\, G_v(t, x; x') \tilde v_0(x')\\
  v(t, x) &= \int \dd x'\, G_v(t, x; x') \ee^{\mu (x-x') } v_0(x') \label{eq:exactsolacgravV1} \\
  \nonumber&=\frac1{2} \int \dd x'\,  \ee^{\mu (x-x') } \del_{ct}  J_0\left(\mu \sqrt{(ct)^2 - (x-x')^2}\right) v_0(x') \\
  \nonumber&+ \frac12 \Big( \ee^{ -\mu ct } v_0(x + ct) +  \ee^{\mu ct } v_0(x -ct) \Big )\\
  \rho(t, x) &= -\frac1{2c} \int \dd x'\,  \ee^{\mu (x-x') } \left( \mu + \del_x \right )  J_0\left(\mu \sqrt{(ct)^2 - (x-x')^2}\right)   v_0(x')\\
  \nonumber& - \frac{1}{2c} \Big( \ee^{ -\mu c t } v_0(x + ct) -  \ee^{\mu c t} v_0(x -ct) \Big ) \label{eq:exactsolacgravV2}
\end{align}
and analogously for $p$. However, it is easier to note that 
\begin{align}
 \del_t (c^2 \rho - p) &= 0
\end{align}
such that
\begin{align}
 p(t, x) = p_0(x) + c^2 \Big(\rho(t, x) - \rho_0(x)\Big)
\end{align}

\end{document}